\documentclass[final,1p,times]{elsarticle}

\usepackage{amsmath}
\usepackage{amssymb}
\usepackage{amsthm}
\allowdisplaybreaks

\usepackage{hyperref}

\usepackage{mathtools}

\usepackage{mismath}

\usepackage{bm} 

\usepackage{graphicx}
\graphicspath{{figs/}} 

\usepackage{tikz}
\usetikzlibrary{math}
\usetikzlibrary{calc}

\usepackage{subcaption}

\usepackage{multirow}

\usepackage{siunitx}
\usepackage{makecell}
\setcellgapes{3pt}

\usepackage{algorithm}

\usepackage{algpseudocode}

\usepackage[capitalize]{cleveref}

\usepackage{easyReview}
\setreviewsoff

\newcommand{\Real}{\mathbb{R}}

\DeclareMathOperator{\Div}{div}

\DeclareMathOperator{\Diag}{diag}
\DeclareMathOperator{\Span}{span}
\DeclareMathOperator{\Cond}{cond}

\DeclareMathOperator{\Tr}{tr}

\DeclarePairedDelimiter{\RoundBrackets}{(}{)}
\DeclarePairedDelimiter{\CurlyBrackets}{\{}{\}}

\usepackage{stmaryrd}
\DeclarePairedDelimiter{\DSquareBrackets}{\llbracket}{\rrbracket}

\newtheorem{theorem}{Theorem}[section]

\newtheorem{lemma}[theorem]{Lemma}

\usepackage{array}
\newcolumntype{L}[1]{>{\raggedright\let\newline\\\arraybackslash\hspace{0pt}}m{#1}}
\newcolumntype{C}[1]{>{\centering\let\newline\\\arraybackslash\hspace{0pt}}m{#1}}
\newcolumntype{R}[1]{>{\raggedleft\let\newline\\\arraybackslash\hspace{0pt}}m{#1}}

\usepackage{algpseudocode}
\algrenewcommand\algorithmiccomment[2][\normalsize]{{#1\hfill\(\triangleright\) #2}}

\usepackage{xcolor}
\definecolor{MyColor1}{HTML}{BFCCB5}
\definecolor{MyColor2}{HTML}{7C96AB}
\definecolor{MyColor3}{HTML}{B7B7B7}
\definecolor{MyColor4}{HTML}{EDC6B1}

\newcommand{\SSSText}[1]{{\scriptscriptstyle\mathup{#1}}}

\journal{arXiv}

\begin{document}

\begin{frontmatter}

    \title{An efficient multiscale multigrid preconditioner for Darcy flow in high-contrast media}


    \author[CUHK]{Changqing Ye}

    \author[EIAS]{Shubin Fu\corref{cor1}}
    \cortext[cor1]{Corresponding author.}
    \ead{sfu@eitech.edu.cn}

    \author[CUHK]{Eric T.~Chung}

    \author[LSEC,AMSS]{Jizu Huang}




    \affiliation[CUHK]{
        organization={Department of Mathematics, The Chinese University of Hong Kong},
        city={Shatin},
        country={Hong Kong SAR}
    }

    \affiliation[EIAS]{
        organization={Eastern Institute for Advanced Study},
        city={Ningbo},
        postcode={315200},
        state={Zhejiang},
        country={PR China}
    }

    \affiliation[LSEC]{
        organization={LSEC, Academy of Mathematics and Systems Science, Chinese Academy of Sciences},
        city={Beijing},
        postcode={100190},
        country={PR China}
    }

    \affiliation[AMSS]{
        organization={School of Mathematical Sciences, University of Chinese Academy of Sciences},
        city={Beijing},
        postcode={100049},
        country={PR China}
    }

    \begin{abstract}

In this paper, we develop a multigrid preconditioner to solve Darcy flow in highly heterogeneous porous media. The key component of the preconditioner is to construct a sequence of nested subspaces $W_{\mathcal{L}}\subset W_{\mathcal{L}-1}\subset\cdots\subset W_1=W_h$. An appropriate spectral problem is defined in the space of $W_{i-1}$, then the eigenfunctions of the spectral problems are utilized to form $W_i$. The preconditioner is applied to solve a positive semidefinite linear system which results from discretizing the Darcy flow equation with the lowest order Raviart-Thomas spaces and adopting a trapezoidal quadrature rule. Theoretical analysis and numerical investigations of this preconditioner will be presented. In particular, we will consider several typical highly heterogeneous permeability fields whose resolutions are up to $1024^3$ and examine the computational performance of the preconditioner in several aspects, such as strong scalability, weak scalability, and robustness against the contrast of the media. We also demonstrate an application of this preconditioner for solving a two-phase flow benchmark problem.

    \end{abstract}




\end{frontmatter}


\section{Introduction}
Simulating fluid flow in porous media is critical for many practical applications, such as reservoir simulation, CO\textsubscript{2} sequestration, nuclear water storage, and underground water contamination. Thanks to the development of reservoir characterization methods, detailed multiscale geocellular models that may contain billions of grid cells become available. Additionally, some simulations, such as CO\textsubscript{2} migration, usually last hundreds of years. As a result, simulating this phenomenon is prohibitively expensive, and it is necessary to perform some model reduction. Generally, there are two categories of model reduction approaches. One is to reduce the model's size by upscaling \cite{Wu2002,Durlofsky1991,Arbogast2013,Chung2023}, which is suitable for applications where small details of the flow can be ignored. Another category is multiscale techniques, including the multiscale finite element method \cite{Hou1997} and its extensions mixed multiscale finite element method \cite{Chen2003,Aarnes2004}, generalized multiscale finite element method (GMsFEM) \cite{Efendiev2013,Chung2015,Chung2023}, multiscale finite volume method \cite{Jenny2003,Hajibeygi2009} and multiscale mortar mixed finite element method \cite{Arbogast2007}. For most of these methods, the basic idea is to solve problems in coarse grids with multiscale basis functions that are constructed by solving carefully defined local problems. These multiscale basis functions are integrated with important information about geological models; hence the accuracy of cheap coarse-grid simulations can be guaranteed.

Although the aforementioned multiscale approaches have been applied successfully to solve a wide range of multiscale problems, multiscale solutions can deteriorate with increasing contrast in permeability and correlation length \cite{Lunati2007,Arbogast2015} in certain scenarios. In cases where highly accurate fine-scale solutions are required, it is necessary to design an efficient solver to solve the original fine-scale problems. Direct solvers such as MUMPS \cite{Amestoy2001} are more applicable for problems with multiple sources since these solvers first factorize target operators and then solve linear systems with forward and back substitutions. However, factorizing large-size matrices requires a huge amount of memory, and, worse, the parallel efficiency of direct solvers is not satisfied. Therefore, for large-scale flow simulations which commonly have fixed sources, adopting preconditioned iterative solvers is preferred. Utilizing multiscale coarse spaces is proven to improve the robustness and efficiency of classical preconditioners, and hence, various novel preconditioners \cite{Graham2007,Sarkis1997,Bjoerstad2002,Wang2014,Calvo2016,Kim2017,Dolean2012,Kim2018,Klawonn2015,Mandel2007,Galvis2010,Galvis2010a,Nataf2011,Klawonn2016,Heinlein2019,Xie2019,Yang2019,AlDaas2021,Xie2014,Bastian2022,Ye2023} have been proposed in past three decades.

Among above mentioned multiscale preconditioners, most of them \cite{Graham2007,Sarkis1997,Bjoerstad2002,Wang2014,Calvo2016,Kim2017,Dolean2012,Kim2018,Klawonn2015,Mandel2007,Galvis2010,Galvis2010a,Nataf2011,Klawonn2016,Heinlein2019,Bastian2022,Kalchev2016} are developed for handling the elliptic problems in second-order formulation, therefore one can not easily obtain mass conservative velocity fields which limits their applications in flow transport problems. In addition, in these works, local spectral problems are constructed in overlapped subdomains, which implies huge communication that deteriorates the efficiency of parallel computing. In \cite{Xie2019,Yang2019,Kalchev2016}, multiscale preconditioners are proposed for solving the saddle system obtained by discretizing the first-order formulation of the elliptic problem, despite their robustness and fast convergence, the degrees of freedom (DoF) of the targeted linear systems are still huge. Moreover, one also needs to solve local spectral problems in overlapped subdomains. To overcome the drawbacks of these multiscale preconditioners, a novel two-grid and two-level multiscale preconditioner is recently developed in \cite{Fu2022} and \cite{Ye2023} where local spectral problems are defined in non-overlapping coarse blocks. In addition, a velocity elimination technique is adopted to reduce the unknowns of the saddle system without hurting the generation of a mass-conservative velocity field. Although this novel efficient preconditioner shows robustness and impressive computational performance, it is difficult to scale up to extreme-scale problems since solving the coarse problems with direct solvers becomes a bottleneck. A potential remedy is to extend this two-grid/two-level preconditioner to the multigrid case which is the main goal of this paper. We will construct a sequence of nested spectral subspaces $W_{\mathcal{L}}\subset W_{\mathcal{L}-1}\subset\cdots\subset W_1=W_h$ and then design a corresponding spectral multigrid preconditioner. The construction of subspace $W_{i}$ relies on solving spectral problems defined in local $W_{i-1}$. We note that similar ideas can be found in \cite{Efendiev2011,Bastian2022,AlDaas2021} which also aim to solve the linear system generated by discretizing the second-order formulation of the high-contrast elliptic problem.

Rich numerical experiments will be presented with several typical highly heterogeneous models. In particular, we examine the strong and weak scalability performance and robustness of our preconditioner against contrast ratios of media. Moreover, the implementation of the proposed preconditioner on a two-phase flow benchmark problem is detailed. We want to mention that multiscale preconditioners are not the unique way to cope with high-contrast flow problems, other types of preconditioners can be found in \cite{Yang2019a,Yang2022,Li2020,Yang2023,Luo2021,Hu2013a,Hu2013,Arbogast2015,Ganis2014,Ganis2012}. Leveraging the convergence theory of inexact two-grid methods, We also provide an analysis of the condition number of the preconditioned operator.

The remainder of the paper is arranged as follows. \Cref{sec:pre} introduces the basic model of subsurface flows, the definition of notation, and fine-scale discretization. \Cref{sec:gms} illustrates a systematic way to construct the coarse space and the preconditioner. An analysis for the proposed preconditioner is presented in \cref{sec:anal}. The numerical results are shown in \cref{sec:num}. The application of a two-phase flow problem is described in \cref{sec:spe10}. Finally, a conclusion is given.

\section{Preliminaries}\label{sec:pre}
\subsection{The model problem and the discretization}
We consider the following elliptic PDE with the unknown pressure field $p$
\begin{equation}\label{eq:ell1}
    -\Div \RoundBrackets*{\mathbb{K}\nabla p}=f \  \text{in} \  \Omega \\
\end{equation}
in a bounded Lipschitz domain $\Omega\in \Real^d$, with $d=2$ or $3$, subject to a no-flow boundary condition
\begin{equation}\label{eq:ell2}
    \mathbb{K}\nabla p\cdot\bm{n}=0 \  \text{on} \ \partial\Omega,
\end{equation}
where $\mathbb{K}$ is the symmetric matrix-valued permeability field that exhibits discontinuity and high heterogeneity, i.e.,
\[
    \frac{\max_{\bm{x}\in \Omega}\lambda_{\mathup{max}}\big(\mathbb{K}(\bm{x})\big)}{\min_{\bm{x}\in \Omega}\lambda_{\mathup{min}}\big(\mathbb{K}(\bm{x})\big)}\gg 1.
\]
In reservoir simulation, the permeability field may not be point-wisely isotropic but commonly be point-wisely orthotropic, that is,
\[
    \mathbb{K}(\bm{x})=\Diag\big(\kappa^x(\bm{x}), \kappa^y(\bm{x})\big)\ \text{or}\ \Diag\big(\kappa^x(\bm{x}), \kappa^y(\bm{x}), \kappa^z(\bm{x})\big),
\]
and we assume the field $\mathbb{K}$ belongs to this class throughout the article. Note that we make a slight abuse of notation for denoting three orthogonal directions by $x$, $y$, and $z$. We emphasize here that the extension of the forthcoming discretization scheme of \cref{eq:ell1} to the full matrix-valued permeability field is possible, while several delicate constructions are necessary \cite{Ingram2010}. The divergence theorem implies that the source function $f$ must satisfy the so-called compatibility condition $\int_{\Omega} f\di \bm{x}=0$, and an additional restriction $\int_{\Omega} p\di \bm{x}=0$ must be imposed to guarantee the uniqueness of the solution. However, if we incorporate well modeling (see \cite{Chen2007,Chen2006}), the source term $f$ will be coupled with pressure and we will revisit this situation in \cref{sec:num}. By introducing the flux variable $\bm{v}=-\mathbb{K}\nabla p$, we can rewrite \cref{eq:ell1} and \cref{eq:ell2} as
\begin{equation}\label{eq:orgional_equation}
    \left\{
    \begin{aligned}
        \mathbb{K}^{-1}\bm{v}+\nabla p & = \bm{0} \ \text{in} \ \Omega,\       &  & \text{(Darcy's law)},                \\
        \Div(\bm{v})                   & = f \ \text{in} \ \Omega, \           &  & \text{(mass conservation)},          \\
        \bm{v}\cdot \bm{n}             & = 0 \ \text{on} \ \partial \Omega, \  &  & \text{(no-flow boundary condition)}.
    \end{aligned}
    \right.
\end{equation}
Note that for Darcy's law in \cref{eq:orgional_equation}, we implicitly normalize the viscosity for simplicity.

The mixed formulation (ref.~\cite{Brezzi1991}) of \cref{eq:orgional_equation} is to seek $(\bm{v}, p)\in \bm{V}\times W$ satisfying the following equations:
\begin{equation}\label{eq:weak}
    \begin{aligned}
        \int_\Omega\mathbb{K}^{-1}\bm{v}\cdot \bm{w} \di \bm{x}-\int_\Omega\Div(\bm{w})p\di \bm{x} & = 0,                         &  & \forall \bm{w}\in \bm{V}, \\
        -\int_\Omega\Div(\bm{v})q \di \bm{x}                                                       & = -\int_\Omega fq\di \bm{x}, &  & \forall q \in W.
    \end{aligned}
\end{equation}
where the space $\bm{V}$ is defined as
\[
    \bm{V}\coloneqq\CurlyBrackets*{\bm{v}\in L^2(\Omega;\Real^d)\mid \Div(\bm{v})\in L^2(\Omega) \ \text{and} \ \bm{v}\cdot\bm{n}=0\;\mbox{on}\;\partial \Omega},
\]
and $W$ is
\[
    W\coloneqq\CurlyBrackets*{q\in L^2(\Omega)\mid \int_{\Omega}q\di \bm{x}=0}.
\]

By restricting the trial and test functions to finite-dimensional subspaces $\bm{V}_h\subset \bm{V}$ and $ W_h\subset W$ associated with a prescribed triangulation $\mathcal{T}_h$ of $\Omega$, the corresponding discrete problem is to find $(\bm{v}_h,p_h)\in \bm{V}_h\times W_h$ such that
\begin{equation}\label{eq:weak2}
    \begin{aligned}
        \int_\Omega\mathbb{K}^{-1}\bm{v}_h\cdot \bm{w}_h\di \bm{x}-\int_\Omega \Div(\bm{w}_h)p_h\di \bm{x} & = 0,                            &  & \forall \bm{w}_h\in \bm{V}_h, \\
        -\int_\Omega\Div(\bm{v}_h)q_h \di \bm{x}                                                           & = -\int_\Omega fq_h \di \bm{x}, &  & \forall q_h \in W_h.
    \end{aligned}
\end{equation}
Several well-known families of mixed finite element spaces $\bm{V}_h \times W_h$ that satisfy the well-known LBB condition can be utilized (see \cite{Auricchio2017}). We choose the lowest-order Raviart-Thomas element space for simplicity, whose definition is detailed in \cite{Brezzi1991,Boffi2013}. Let $\CurlyBrackets*{\bm{\phi}_e}_{e\in \mathcal{E}^0_h}$ and $\CurlyBrackets*{q_\tau}_{\tau\in \mathcal{T}_h}$ be the bases of $\bm{V}_h$ and $W_h$, respectively, where $\mathcal{E}^0_h$ is the set of internal edges/faces, and we drop $h$ of these bases for brevity. Hence, the finite element solution $(\bm{v}_h, p_h)$ can be expressed as
\[
    \bm{v}_h=\sum_{e\in\mathcal{E}^0_h}v_e\bm{\phi}_e \ \text{and} \ p_h=\sum_{\tau\in\mathcal{T}_h}p_\tau q_\tau.
\]
Then, it is easy to obtain the following algebraic representation of \cref{eq:weak2}
\begin{equation}\label{eq:fine_system}
    \begin{pmatrix}
        \mathsf{M}  & -\mathsf{B}^\intercal \\
        -\mathsf{B} & \mathsf{0}            \\
    \end{pmatrix}
    \begin{pmatrix}
        \mathsf{v} \\ \mathsf{p}
    \end{pmatrix}
    =\begin{pmatrix}
        \mathsf{0} \\
        -\mathsf{f}
    \end{pmatrix},
\end{equation}
where $\mathsf{M}$ is a symmetric, positive definite matrix, $\mathsf{v}$, $\mathsf{p}$ and $\mathsf{f}$ are column vectors. 

\subsection{The velocity elimination technique}\label{subsec:vet}
Preconditioning the saddle point system \cref{eq:fine_system} directly is not easy, and we may involve several sophisticated field-splitting procedures. Fortunately, if applying the trapezoidal quadrature rule to compute the integration $\int_{\Omega}\mathbb{K}^{-1}\bm{v}_h\cdot\bm{w}_h\di \bm{x}$ on each element, we will obtain a diagonal velocity mass matrix $\mathsf{M}_\mathup{t}$, and hence \cref{eq:fine_system} could be replaced by
\[
    \begin{pmatrix}
        \mathsf{M}_\mathup{t} & -\mathsf{B}^\intercal \\
        -\mathsf{B}           & \mathsf{0}            \\
    \end{pmatrix}
    \begin{pmatrix}
        \mathsf{v} \\ \mathsf{p}
    \end{pmatrix}
    =\begin{pmatrix}
        \mathsf{0} \\
        -\mathsf{f}
    \end{pmatrix}.
\]
Since $\mathsf{M}_\mathup{t}$ could be inverted straightforwardly and the unknown $\mathsf{v}$ in \cref{eq:fine_system} can be fully eliminated, we boil down to solve a linear system of $\mathsf{p}$ as
\begin{equation}\label{eq:p}
    \mathsf{A}\mathsf{p}=\mathsf{f},
\end{equation}
where $\mathsf{A}\coloneqq\mathsf{B}\mathsf{M}_\mathup{t}^{-1}\mathsf{B}^\intercal$ is a symmetric and semi-positive definite matrix. This technique is known as velocity elimination \cite{Russell1983,Arbogast1997,Chen2020}, which could also be interpreted as a cell-centered finite difference scheme on Cartesian meshes. Once $\mathsf{p}$ is solved, the velocity field could be recovered by $\mathsf{v}=\mathsf{M}_\mathup{t}^{-1}\mathsf{B}^\intercal \mathsf{p}$.

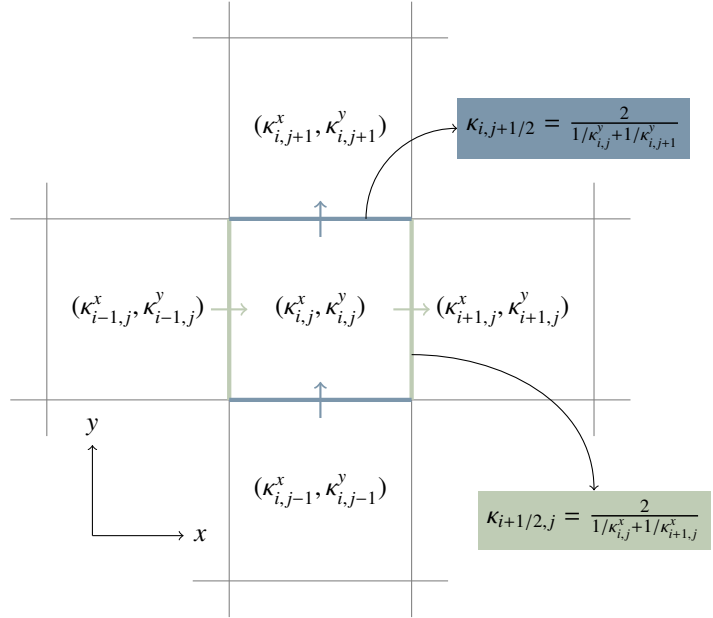
\begin{figure}[!ht]
    \centering
    \begin{tikzpicture}[scale=1.2]
        \draw[gray] (1.6, 0.0) -- (4.4, 0.0);
        \draw[gray] (-0.4, 2.0) -- (6.4, 2.0);
        \draw[gray] (-0.4, 4.0) -- (6.4, 4.0);
        \draw[gray] (1.6, 6.0) -- (4.4, 6.0);
        \draw[gray] (0.0, 1.6) -- (0.0, 4.4);
        \draw[gray] (2.0, -0.4) -- (2.0, 6.4);
        \draw[gray] (4.0, -0.4) -- (4.0, 6.4);
        \draw[gray] (6.0, 1.6) -- (6.0, 4.4);

        \node at (3.0, 3.0) {$(\kappa^x_{i,j}, \kappa^y_{i,j})$};
        \node at (1.0, 3.0) {$(\kappa^x_{i-1,j}, \kappa^y_{i-1,j})$};
        \node at (5.0, 3.0) {$(\kappa^x_{i+1,j}, \kappa^y_{i+1,j})$};
        \node at (3.0, 1.0) {$(\kappa^x_{i,j-1}, \kappa^y_{i,j-1})$};
        \node at (3.0, 5.0) {$(\kappa^x_{i,j+1}, \kappa^y_{i,j+1})$};

        \draw[->, thick, MyColor1] (1.8, 3.0) -- ++(0.4, 0.0);
        \draw[->, thick, MyColor1] (3.8, 3.0) -- ++(0.4, 0.0);
        \draw[->, thick, MyColor2] (3.0, 1.8) -- ++(0.0, 0.4);
        \draw[->, thick, MyColor2] (3.0, 3.8) -- ++(0.0, 0.4);

        \draw[ultra thick, MyColor1] (2.0, 2.0) -- ++(0.0, 2.0);
        \draw[ultra thick, MyColor1] (4.0, 2.0) -- ++(0.0, 2.0);
        \draw[ultra thick, MyColor2] (2.0, 2.0) -- ++(2.0, 0.0);
        \draw[ultra thick, MyColor2] (2.0, 4.0) -- ++(2.0, 0.0);

        \draw[->, thin] (0.5, 0.5) -- ++(1.0, 0.0) node[right] {$x$};
        \draw[->, thin] (0.5, 0.5) -- ++(0.0, 1.0) node[above] {$y$};

        \draw[->] (3.5, 4.0) to [out=90, in=180] (4.5, 5.0) node[right,fill=MyColor2] {$\kappa_{i,j+1/2}=\frac{2}{1/\kappa_{i,j}^y+1/\kappa_{i,j+1}^y}$};
        \draw[->] (4.0, 2.5) to [out=0, in=90] (6.0, 1.0) node[below,fill=MyColor1] {$\kappa_{i+1/2,j}=\frac{2}{1/\kappa_{i,j}^x+1/\kappa_{i+1,j}^x}$};
    \end{tikzpicture}
    \caption{An illustration of the expression of $\kappa_e$. Depending on the direction, $\kappa_{i,j+1/2}$ and $\kappa_{i+1/2,j}$ utilize harmonic averages of different component of the permeability field.}\label{fig:kappa edge}
\end{figure}

Because the analysis result in \cref{sec:anal} is based on the aforementioned scheme, we elaborate on several details of the discretization here. We limit our discussion to a 2D rectangular domain and structured or Cartesian grids. Let the size of each element be $h^x\times h^y$, the variational form corresponding to the linear system \cref{eq:p} is
\begin{equation}\label{eq:varia vet}
    \sum_{e\in\mathcal{E}_h^0} \kappa_e \DSquareBrackets{p_h}_e\DSquareBrackets{q_h}_e\frac{\abs{e}^2}{h_xh_y} = \int_\Omega f q_h \di \bm{x}, \ \forall q_h \in W_h,
\end{equation}
where $\abs{e}$ is the length of the edge $e$ with taking the value of $h^x$ or $h^y$. The two notations $\kappa_e$ and $\DSquareBrackets{\cdot}_e$ in \cref{eq:varia vet} deserve to be addressed: for each internal edge $e$, $\DSquareBrackets{p_h}_e$ is the jump of $p_h$ across $e$ towards the positive $x$- or $y$-direction; the value of $\kappa_e$ is illustrated in \cref{fig:kappa edge}, and the harmonic average stands out as an important property. If we identify $f$ as a function in $W_h$, the variational form \cref{eq:varia vet} enjoys an pure algebraic representation as
\begin{equation}\label{eq:discrete varia}
    \begin{aligned}
         & \sum_{i,j}\frac{\kappa_{i+1/2,j}}{(h^x)^2}\RoundBrackets*{p_{i+1,j}-p_{i,j}}\RoundBrackets*{q_{i+1,j}-q_{i,j}}                        \\
         & \qquad+\frac{\kappa_{i,j+1/2}}{(h^y)^2}\RoundBrackets*{p_{i,j+1}-p_{i,j}}\RoundBrackets*{q_{i,j+1}-q_{i,j}}=\sum_{i,j}f_{i,j}q_{i,j},
    \end{aligned}
\end{equation}
and note here we cancel out $h^xh^y$ the area of a fine element on both sides of \cref{eq:varia vet}. The conservation of mass on the fine element could be derived from \cref{eq:discrete varia} as
\[
    h^y\RoundBrackets*{v_{i+1/2,j}-v_{i-1/2,j}}+h^x\RoundBrackets*{v_{i,j+1/2}-v_{i,j-1/2}}=h^xh^yf_{i,j}
\]
with
\begin{equation}\label{eq:velocity}
    v_{i+1/2,j}=-\kappa_{i+1/2,j}\frac{p_{i+1,j}-p_{i,j}}{h^x}\ \text{and} \ v_{i,j+1/2}=-\kappa_{i,j+1/2}\frac{p_{i,j+1}-p_{i,j}}{h^y}.
\end{equation}
This scheme shares similarities with finite-volume methods by treating each fine element as a control volume and is termed the two-point flux approximation scheme in the literature \cite{Barth2017}. By normalizing $\mathbb{K}$ as
\[
    \widetilde{\mathbb{K}}=\Diag\RoundBrackets*{\widetilde{\kappa}^x, \widetilde{\kappa}^y\big)\coloneqq \Diag\big({\kappa^x}/{(h^x)^2}, {\kappa^y}/{(h^y)^2}},
\]
we can see that
\[
    \frac{\kappa_{i+1/2,j}}{(h^x)^2}=\widetilde{\kappa}_{i+1/2,j}=\frac{2}{1/\widetilde{\kappa}_{i,j}^x+1/\widetilde{\kappa}_{i+1,j}^x}\ \text{and} \ \frac{\kappa_{i,j+1/2}}{(h^y)^2}=\widetilde{\kappa}_{i,j+1/2}=\frac{2}{1/\widetilde{\kappa}_{i,j}^y+1/\widetilde{\kappa}_{i,j+1}^y},
\]
which frees \cref{eq:discrete varia} of dimensions of fine elements. This transformation is physically reasonable because the permeability $\mathbb{K}$ is expressed in units of $[\mathup{Length}]^2$ while $\widetilde{\mathbb{K}}$ is dimensionless. We will utilize this notation in the construction of our preconditioner later.

The main goal of this paper is to develop a spectral multigrid preconditioner for solving \cref{eq:p}.

\section{Multiscale coarse space based multigrid preconditioner}\label{sec:gms}
To present the construction of coarse spaces, we first explain the decomposition of the domain $\Omega$. We assume that $\Omega$ is meshed by $\mathcal{L}$ layers of hierarchical grids that are denoted by $\mathcal{T}^{(l)}$ with $l\in \CurlyBrackets*{0,\dots,\mathcal{L}-1}$. The hierarchy implies that each $\mathcal{T}^{(i-1)}$ is a nested refinement of $\mathcal{T}^{(i)}$ for $i=1,\dots,\mathcal{L}-1$. Recalling that the finest mesh is $\mathcal{T}_h=\mathcal{T}^{(0)}$ and the resolution of $\mathbb{K}$ is tied to $\mathcal{T}_h$, We shall treat $h$ as a given physical parameter rather than a variable index in convergence theories of finite element methods.

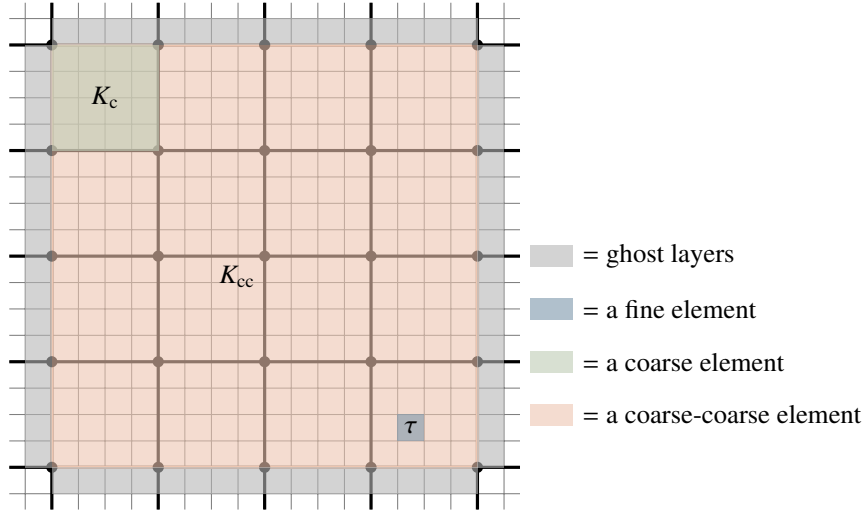
\begin{figure}[!ht]
    \centering
    \begin{tikzpicture}[scale=1.4]
        \draw[step=0.25, gray, thin] (-0.4, -0.4) grid (4.4, 4.4);
        \draw[step=1.0, black, very thick] (-0.4, -0.4) grid (4.4, 4.4);
        \draw[MyColor4, very thick] (0.0, 0.0) rectangle (4.0, 4.0);
        \foreach \x in {0,...,4}
        \foreach \y in {0,...,4}{
                \fill (1.0 * \x, 1.0 * \y) circle (1.5pt);
            }

        \fill[MyColor4, opacity=0.6] (0.0, 0.0) rectangle (4.0, 4.0);
        \fill[MyColor1, opacity=0.6] (0.0, 3.0) rectangle (1.0, 4.0);
        \fill[MyColor2, opacity=0.6] (3.25, 0.25) rectangle (3.5, 0.5);
        \fill[MyColor3, opacity=0.6] (-0.25, 0.0) rectangle (0.0, 4.0);
        \fill[MyColor3, opacity=0.6] (0.0, 0.0) rectangle (4.0, -0.25);
        \fill[MyColor3, opacity=0.6] (0.0, 4.0) rectangle (4.0, 4.25);
        \fill[MyColor3, opacity=0.6] (4.25, 0.0) rectangle (4.0, 4.0);

        \node at (0.5, 3.5) {$K_\mathup{c}$};
        \node[below left] at (2.0, 2.0) {$K_\mathup{cc}$};
        \node at (3.375, 0.375) {$\tau$};

        \fill[MyColor3, opacity=0.6] (4.5, 1.9) rectangle (4.9, 2.1);
        \node[right] at (4.9, 2.0) {=\ ghost layers};

        \fill[MyColor2, opacity=0.6] (4.5, 1.4) rectangle (4.9, 1.6);
        \node[right] at (4.9, 1.5) {=\ a fine element};

        \fill[MyColor1, opacity=0.6] (4.5, 0.9) rectangle (4.9, 1.1);
        \node[right] at (4.9, 1.0) {=\ a coarse element};

        \fill[MyColor4, opacity=0.6] (4.5, 0.4) rectangle (4.9, 0.6);
        \node[right] at (4.9, 0.5) {=\ a coarse-coarse element};
    \end{tikzpicture}
    \caption{An illustration of hierarchical meshes, a fine element $\tau$, a coarse element $K_\mathup{c}$, a coarse-coarse element $K_\mathup{cc}$ and ghost layers. In MPI-implementations, each $K_\mathup{cc}$ is assigned to a unique process.}
    \label{fig:grid}
\end{figure}

In this paper, we focus on a three-grid setting, i.e., $\mathcal{L}=3$. For simplicity of notations, we denote by $\mathcal{T}_\mathup{c}=\mathcal{T}^{(1)}$ the coarse mesh and $\mathcal{T}_\mathup{cc}=\mathcal{T}^{(2)}$ the coarse-coarse mesh. This setting is also based on the MPI parallel computing architecture: each MPI process handles a coarse-coarse element in $\mathcal{T}_\mathup{cc}$ with ghost layers of the width of one fine element, where ghost layers specialize in inter-process communications; each coarse-coarse element is further partitioned into coarse elements, and there are no inter-process communications required between them. \Cref{fig:grid} is an illustration of hierarchical meshes, where a fine element $\tau$ in $\mathcal{T}_h$, a coarse element $K_\mathup{c}$ in $\mathcal{T}_\mathup{c}$, a coarse-coarse element $K_\mathup{cc}$ in $\mathcal{T}_\mathup{cc}$ and ghost layers are visualized in different colors. We denote by $m_\mathup{c}$ the total number of coarse elements, $m_\mathup{cc}$ the total number of coarse-coarse elements, and $n$ the total number of fine elements (also DoF).

A typical three-grid iteration is demonstrated in \cref{alg:three grid}. Although, under several assumptions, three-grid or multigrid methods are proven to be convergent (see, e.g., \cite{Notay2007,Xu2022}), we, however, exploit them in preconditioned iterative solvers as an acceleration technique (ref.~\cite{Saad2003}). The matrices $\mathsf{M}$ and $\mathsf{M}_\mathup{c}$ are the smoothers on the fine and coarse grids, respectively. The two restriction matrices $\mathsf{R}_\mathup{c}$ and $\mathsf{R}_\mathup{cc}$ are crucial ingredients of our preconditioner, and the transposes $\mathsf{R}_\mathup{c}^\intercal$ and $\mathsf{R}_\mathup{cc}^\intercal$ are referred to as prolongation matrices conventionally. Moreover, we have the relations $\mathsf{A}_\mathup{c}=\mathsf{R}_\mathup{c}\mathsf{A}\mathsf{R}_\mathup{c}^\intercal$ and $\mathsf{A}_\mathup{cc}=\mathsf{R}_\mathup{cc}\mathsf{A}_\mathup{c}\mathsf{R}_\mathup{cc}^\intercal=(\mathsf{R}_\mathup{cc}\mathsf{R}_\mathup{c})\mathsf{A}(\mathsf{R}_\mathup{cc}\mathsf{R}_\mathup{c})^\intercal$, where $\mathsf{A}_\mathup{cc}$ could be singular, and hence the pseudoinverse $\mathsf{A}_\mathup{cc}^\dagger$ is involved in \cref{alg:three grid}. The three-grid method differentiates itself from the two-grid version by the map $\mathsf{r}_\mathup{c}\mapsto \mathsf{e}_\mathup{c}$ in \cref{alg:three grid}, where the latter considers an exact inversion as $\mathsf{e}_\mathup{c}=\mathsf{A}_\mathup{c}^\dagger\mathsf{r}_\mathup{c}$. In our previous effort \cite{Ye2023} on a two-level overlapping preconditioner, we emphasize that the coarse solver is a bottleneck for the computing performance to be scaled up. The aforementioned three-grid method, in some sense, could be viewed as an inexact two-grid method by defining an operator $\mathsf{B}_\mathup{c}^{-1}$ in the following form
\begin{equation}\label{eq:B_c}
    \begin{aligned}
        \mathsf{e}_\mathup{c} & = \mathsf{B}_\mathup{c}^{-1}\mathsf{r}_\mathup{c}                                                                                                                                                                                                                                                                                                                                                                                                \\
                              & \coloneqq \Big\{ \mathsf{M}_\mathup{c}^{-\intercal}+\mathsf{M}_\mathup{c}^{-1}-\mathsf{M}_\mathup{c}^{-\intercal}\mathsf{A}_\mathup{c}\mathsf{M}_\mathup{c}^{-1}+\RoundBrackets*{\mathsf{I}-\mathsf{M}_\mathup{c}^{-\intercal}\mathsf{A}_\mathup{c}}\mathsf{R}_\mathup{cc}^\intercal\mathsf{A}_\mathup{cc}^{\dagger}\mathsf{R}_\mathup{cc}\RoundBrackets*{\mathsf{I}-\mathsf{A}_\mathup{c}\mathsf{M}_\mathup{c}^{-1}}\Big\}\mathsf{r}_\mathup{c}
    \end{aligned}
\end{equation}
as an approximation to $\mathsf{A}_\mathup{c}^\dagger$.

\begin{algorithm}
    \caption{An iteration in three-grid method (preconditioner).}\label{alg:three grid}
    \begin{algorithmic}[1]
        \Require{The operators---$\mathsf{A}$, $\mathsf{A}_\mathup{c}$, $\mathsf{M}^{-1}$, $\mathsf{M}^{-1}_\mathup{c}$, $\mathsf{R}_\mathup{c}$ and $\mathsf{R}_\mathup{cc}$; the right-hand vector---$\mathsf{f}$; an initial guess---$\mathsf{u}^{(0)}$}
        \State{Presmoothing: $\mathsf{u}^{(1)}\leftarrow \mathsf{u}^{(0)}+\mathsf{M}^{-1}(\mathsf{f}-\mathsf{A}\mathsf{u}^{(0)})$}  \Comment[\footnotesize]{$\mathsf{M}$ is the smoother on the fine grid}
        \State{Restriction: $\mathsf{r}_\mathup{c}\leftarrow \mathsf{R}_\mathup{c}(\mathsf{f}-\mathsf{A}\mathsf{u}^{(1)})$}  \Comment[\footnotesize]{Restrict the residual onto the coarse space}
        \State{\quad Presmoothing: $\mathsf{u}_c^{(0)}\leftarrow \mathsf{M}_\mathup{c}^{-1}\mathsf{r}_\mathup{c}$}  \Comment[\footnotesize]{$\mathsf{M}_\mathup{c}$ is the smoother on the coarse grid}
        \State{\quad Restriction: $\mathsf{r}_\mathup{cc}\leftarrow \mathsf{R}_\mathup{cc}(\mathsf{r}_\mathup{c}-\mathsf{A}_\mathup{c}\mathsf{u}_\mathup{c}^{(0)})$} \Comment[\footnotesize]{Restrict the residual onto the coarse-coarse space}
        \State{\quad\quad Correction: $\mathsf{e}_\mathup{cc}\leftarrow \mathsf{A}_\mathup{cc}^{\dagger}\mathsf{r}_\mathup{cc}$} \Comment[\footnotesize]{Involve a direct solver on the coarse-coarse grid}
        \State{\quad Prolongation: $\mathsf{u}_\mathup{c}^{(1)}\leftarrow \mathsf{u}_\mathup{c}^{(0)}+\mathsf{R}_\mathup{cc}^\intercal\mathsf{e}_\mathup{cc}$ \Comment[\footnotesize]{Project back into the coarse space}}
        \State{\quad Postsmoothing: $\mathsf{e}_\mathup{c}\leftarrow \mathsf{u}_\mathup{c}^{(1)}+\mathsf{M}_\mathsf{c}^{-\intercal}(\mathsf{r}_\mathup{c}-\mathsf{A}_\mathup{c}\mathsf{u}_\mathup{c}^{(1)})$} \Comment[\footnotesize]{Build the map $\mathsf{r}_\mathup{c}\mapsto \mathsf{e}_\mathup{c}$}
        \State{Prolongation: $\mathsf{u}^{(2)}\leftarrow \mathsf{u}^{(1)}+\mathsf{R}_\mathup{c}^\intercal \mathsf{e}_\mathup{c}$}\Comment[\footnotesize]{Project back into the fine space}
        \State{Postsmoothing: $\mathsf{u}_\star\leftarrow \mathsf{u}^{(2)}+\mathsf{M}^{-\intercal}(\mathsf{f}-\mathsf{A}\mathsf{u}^{(2)})$} \Comment[\footnotesize]{The end of the iteration}
    \end{algorithmic}
\end{algorithm}

We then detail the constructions of $\mathsf{R}_\mathup{c}$ and $\mathsf{R}_\mathup{cc}$. For any coarse element $K_\mathup{c}^i\in \mathcal{T}_\mathup{c}$, we define $W_h(K_\mathup{c}^i)$ by restricting $W_h$ on $K_\mathup{c}^i$ and solve a local spectral problem:
\begin{equation}\label{eq:spepb}
    \sum_{e\in \mathcal{E}_h^0(K_\mathup{c}^i)}\kappa_e\DSquareBrackets{\Phi_h}_e\DSquareBrackets{q_h}_e\frac{\abs{e}^2}{h_xh_y} = \lambda \int_{K_\mathup{c}^i} \RoundBrackets*{\Tr \widetilde{\mathbb{K}}} \Phi_h q_h \di \bm{x}, \ \forall q_h \in W_h(K_\mathup{c}^i),
\end{equation}
where $\mathcal{E}_h^0(K_\mathup{c}^i)$ is the set of all internal edges in $K_\mathup{c}^i$, $\Phi_h$ is an eigenvector corresponding to the eigenvalue $\lambda$ and $\widetilde{\mathbb{K}}$ is the normalized permeability field. Note that the eigenvalues here are all dimensionless. Compared to \cite{Ye2023}, we modify the right-hand bilinear form because the permeability field is now assumed to be orthotropic rather than isotropic. After solving the spectral problem \cref{eq:spepb}, we construct the local coarse space $W_\mathup{c}(K_\mathup{c}^i)$ as $\Span\CurlyBrackets{\Phi_{h}^{i,j}\mid j=0,\dots,l_\mathup{c}^i-1}$, where $\CurlyBrackets{\Phi_{h}^{i,j}}_{j=0}^{l_\mathup{c}^{i}-1}$ are the eigenvectors associated with the smallest eigenvalues $L^i_\mathup{c}$. By identifying $W_\mathup{c}(K_\mathup{c}^i)$ as a linear subspace of $W_h$, the global coarse space $W_\mathup{c}\subset W_h$ is formed by $ W_\mathup{c}=W_\mathup{c}(K_\mathup{c}^1)\oplus \dots \oplus W_\mathup{c}(K_\mathup{c}^{m_\mathup{c}})$, and $n_\mathup{c}$---the dimension of $W_\mathup{c}$---is determined by $n_\mathup{c}=l_\mathup{c}^1+\dots+l_\mathup{c}^{m_\mathup{c}}$. Therefore, we establish a $\mathsf{R}_\mathup{c}\in \Real^{n_\mathup{c}\times n}$ as a restriction matrix from the fine grid to the coarse grid, where each row of it is the algebraic representation of an eigenvector of \cref{eq:spepb} in $\Real^n$.

To perform a further dimension reduction on $W_\mathup{c}$, we propose another spectral problem on each coarse-coarse element $K_\mathup{cc}^i \in \mathcal{T}_\mathup{cc}$:
\begin{equation}\label{eq:spectral cc}
    \sum_{e\in \mathcal{E}_h^0(K_\mathup{cc}^i)}\kappa_e\DSquareBrackets{\varPsi_\mathup{c}}_e\DSquareBrackets{\varTheta_\mathup{c}}_e\frac{\abs{e}^2}{h_xh_y} = \lambda \int_{K_\mathup{cc}^i} \RoundBrackets*{\Tr \widetilde{\mathbb{K}}} \varPsi_\mathup{c} \varTheta_\mathup{c} \di \bm{x}, \ \forall \varTheta_\mathup{c} \in W_\mathup{c}(K_\mathup{cc}^i),
\end{equation}
where $\mathcal{E}_h^0(K_\mathup{cc}^i)$ is the set of all internal edges in $K_\mathup{cc}^i$, $W_\mathup{c}(K_\mathup{cc}^i)$ is the restriction of $W_\mathup{c}$ on $K_\mathup{cc}^i$. It should be addressed that \cref{eq:spectral cc} is defined in the linear space $W_\mathup{c}(K_\mathup{cc}^i)$ rather than $W_h(K_\mathup{cc}^i)$, and this fact may be favorable in two aspects: 1) solving \cref{eq:spectral cc} leads to a smaller algebraic system and should be much easier compared with the latter, and 2) the eigenvectors of \cref{eq:spectral cc} compress the coarse space, while the latter one is irrelevant with $W_\mathup{c}$. Similarly, the local coarse-coarse space $W_\mathup{cc}(K_\mathup{cc}^i)$ could be created from spanning $\CurlyBrackets*{\varPsi_\mathup{c}^{i,j}}_{j=0}^{l_\mathup{cc}^i-1}$ which are eigenvectors corresponding to the smallest $l_\mathup{cc}^i$ eigenvalues. Furthermore, by identifying $W_\mathup{cc}(K_\mathup{cc}^i)$ as a linear subspace of $W_\mathup{c}$, we can build the coarse-coarse space $W_\mathup{cc}$ as $W_\mathup{cc}(K_\mathup{cc}^1)\oplus\dots\oplus W_\mathup{cc}(K_\mathup{cc}^{m_\mathup{cc}})$, whose dimension $n_\mathup{cc}$ equals $l_\mathup{cc}^1+\dots+l_\mathup{cc}^{m_\mathup{cc}}$, and thus we achieve a chain of inclusions $W_\mathup{cc}\subset W_\mathup{c} \subset W_h$. Back to $\mathsf{R}_\mathup{cc}$, we let each row of $\mathsf{R}_\mathup{cc}\mathsf{R}_\mathup{c}$ be the algebraic representation eigenvector $\varPsi_c$ of \cref{eq:spectral cc} in $\Real^n$, and  $\mathsf{R}_\mathup{cc}$ has a shape of $n_\mathup{cc}\times n_\mathup{c}$.

Note that the smallest eigenvalue of \cref{eq:spepb,eq:spectral cc} is always zero and the corresponding eigenfunctions take a constant value. This property implies the kernel space of the left-hand operator in \cref{eq:varia vet} is contained in $W_\mathup{c}$ and also $W_\mathup{cc}$. Moreover, considering the Raviart-Thomas spaces on $\mathcal{T}_\mathup{c}$ and $\mathcal{T}_\mathup{cc}$, the pressure parts take a constant value on each element, we can hence interpret $W_\mathup{c}$ and $W_\mathup{cc}$ as enhancements to classic geometric coarse spaces that are based on nested grids.

Another important fact that could be observed from \cref{eq:spepb,eq:spectral cc} is that the right-hand operators in algebraic forms are all diagonal matrices. Moreover, if eigenvectors obtained in \cref{eq:spepb} are normalized w.r.t.~the right-hand bilinear form and are chosen as the bases of $W_\mathup{c}$, the right-hand operator in the algebraic form of \cref{eq:spectral cc} is an identity matrix, which reduces a generalized eigenvalue problem to a standard one. Those properties are thanks to that the right-hand bilinear forms of \cref{eq:spepb,eq:spectral cc} are essentially a weighted $L^2$ inner-product and our coarse elements do not overlap with each other. Constructing coarse spaces via solving spectral problems is criticized as time-consuming (see \cite{Xu2017}), while the proposed spectral problems with a simple right-hand matrix could be expected to alleviate the difficulty.

\section{Analysis}\label{sec:anal}
To avoid working on singular systems, we suppose that a part of $\partial \Omega$ in \cref{eq:orgional_equation} is imposed with a Dirichlet boundary condition, which yields a positive definite matrix $\mathsf{A}$ in \cref{eq:p}. We admit there may exist several subtle modifications for the following theoretical augments to be applied on a positive semidefinite matrix (cf.~\cite{Ye2023}). Nevertheless, the numerical experiments that will be presented in the next section support the applicability of the proposed preconditioner on the original problem \cref{eq:orgional_equation}.

\Cref{alg:three grid} determines an inexact two-grid preconditioner $\mathsf{P}_\SSSText{ITG}$ as
\begin{equation} \label{eq:B_ITG}
    \mathsf{P}_\SSSText{ITG}^{-1}\coloneqq\mathsf{M}^{-\intercal}+\mathsf{M}^{-1}-\mathsf{M}^{-\intercal}\mathsf{A}\mathsf{M}^{-1}+\RoundBrackets*{\mathsf{I}-\mathsf{M}^{-\intercal}\mathsf{A}}\mathsf{R}_\mathup{c}^\intercal\mathsf{B}_\mathup{c}^{-1}\mathsf{R}_\mathup{c}\RoundBrackets*{\mathsf{I}-\mathsf{A}\mathsf{M}^{-1}},
\end{equation}
where the definition of $\mathsf{B}_\mathup{c}^{-1}$ could be referred to \cref{eq:B_c}. Accordingly, the exact preconditioner $\mathsf{P}_\SSSText{TG}$ could be defined as
\begin{equation}\label{eq:B_TG}
    \mathsf{P}_\SSSText{TG}^{-1}\coloneqq\mathsf{M}^{-\intercal}+\mathsf{M}^{-1}-\mathsf{M}^{-\intercal}\mathsf{A}\mathsf{M}^{-1}+\RoundBrackets*{\mathsf{I}-\mathsf{M}^{-\intercal}\mathsf{A}}\mathsf{R}_\mathup{c}^\intercal\mathsf{A}_\mathup{c}^{-1}\mathsf{R}_\mathup{c}\RoundBrackets*{\mathsf{I}-\mathsf{A}\mathsf{M}^{-1}},
\end{equation}
and note that $\mathsf{B}_\mathup{c}^{-1}$ in \cref{eq:B_ITG} is replaced by $\mathsf{A}_\mathup{c}^{-1}$ here. Because $\mathsf{A}$ is a positive definite matrix now, we can safely say $\mathsf{A}_\mathup{c}$ and $\mathsf{A}_\mathup{cc}$ are invertible. The objective of this section is to provide an upper bound of the condition number of the preconditioned system, that is,
\[
    \Cond (\mathsf{P}_\SSSText{ITG}^{-1}\mathsf{A})\coloneqq \frac{\lambda_{\mathup{max}}(\mathsf{P}_\SSSText{ITG}^{-1}\mathsf{A})}{\lambda_{\mathup{min}}(\mathsf{P}_\SSSText{ITG}^{-1}\mathsf{A})} \leq C.
\]
Since the inexact preconditioner is a perturbation of the exact preconditioner, the results in \cite{Notay2007} suggest estimating
\[
    \Cond (\mathsf{P}_\SSSText{TG}^{-1}\mathsf{A})\coloneqq \frac{\lambda_{\mathup{max}}(\mathsf{P}_\SSSText{TG}^{-1}\mathsf{A})}{\lambda_{\mathup{min}}(\mathsf{P}_\SSSText{TG}^{-1}\mathsf{A})}
\]
first.

According to \cite{Xu2022}, we can see that $\lambda_{\mathup{max}}(\mathsf{P}_\SSSText{TG}^{-1}\mathsf{A})=1$. For a positive define matrix $\mathsf{Z}$, we take $\norm{\mathsf{v}}_\mathsf{Z}=\sqrt{\mathsf{v}\cdot \mathsf{Z}\mathsf{v}}$ as the energy norm induced by $\mathsf{Z}$ for a compatible vector $\mathsf{v}$. We usually assume that the smoother $\mathsf{M}$ is properly chosen such that $\mathsf{M}+\mathsf{M}^\intercal-\mathsf{A}$ is a positive definite matrix. According to \cite{Xu2017}, a direct result of this assumption is that two-grid iterations can converge without coarse space corrections, although the convergence rate may be poor. Thanks to the XZ-identity \cite{Xu2017}, we could obtain that
\begin{equation}\label{eq:XZ identity}
    \frac{1}{\lambda_{\mathup{min}}(\mathsf{P}_\SSSText{TG}^{-1}\mathsf{A})}=\max_{\mathsf{v}\in \Real^n\setminus\{\mathsf{0}\}}\min_{\mathsf{v}_\mathup{c}\in \Real^{n_\mathup{c}}\setminus\{\mathsf{0}\}}\frac{\norm{\mathsf{v}-\mathsf{R}^\intercal_\mathup{c}\mathsf{v}_\mathup{c}}_{\widetilde{\mathsf{M}}}^2}{\norm{\mathsf{v}}_\mathsf{A}^2},
\end{equation}
where $\widetilde{\mathsf{M}}=\mathsf{M}^\intercal(\mathsf{M}+\mathsf{M}^\intercal-\mathsf{A})^{-1}\mathsf{M}$. The special max-min form of \cref{eq:XZ identity} inspires an existence of the optimal coarse space, which is an eigenspace of the generalized eigenvalue problem $\mathsf{A}\mathsf{v}=\lambda \widetilde{\mathsf{M}}\mathsf{v}$. Nevertheless, due to the facts that $\widetilde{\mathsf{M}}$ requires an inversion of $\mathsf{M}+\mathsf{M}^\intercal-\mathsf{A}$ and the eigenvalue problem needs to be solved in $\Real^n$,  this optimal coarse space is highly unobtainable.

For a coarse element $K_\mathup{c}^i \in \mathcal{T}_\mathup{c}$, we take $\mathsf{A}_i$ as the matrix built from the left-hand bilinear form of \cref{eq:spepb} and $\mathsf{S}^i$ from the right-hand bilinear form, where $\mathsf{S}_i$ is a diagonal matrix. We have seen that the coarse grid $\mathcal{T}_\mathup{c}$ consists of non-overlapping coarse elements. It is natural to introduce binary matrices $\CurlyBrackets*{\mathsf{E}_1,\dots,\mathsf{E}_{m_\mathup{c}}}$ with the relation $\mathsf{I}=\mathsf{E}_1^\intercal\mathsf{E}_1+\dots+\mathsf{E}_{m_\mathup{c}}^\intercal\mathsf{E}_{m_\mathup{c}}$, where the effect of $\mathsf{E}_i$ is restricting the DoF of the domain $\Omega$ onto the coarse element $K_\mathup{c}^i$. We have the following relation:
\[
    \sum_{i=1}^{m_\mathup{c}} \mathsf{E}_i^\intercal \mathsf{A}_i \mathsf{E}_i\lesssim \mathsf{A},
\]
where $\mathsf{A}'\lesssim \mathsf{A}$ means that $\mathsf{A}-\mathsf{A}'$ is positive semi-definite. This inequality is from taking a summation of the left-hand parts of \cref{eq:spepb} w.r.t.~all coarse elements, the edges in boundaries of coarse elements will not be counted. Therefore, we can derive that
\[
    \norm{\mathsf{v}}_\mathsf{A}^2\geq \sum_{i=1}^{m_\mathup{c}} \RoundBrackets*{\mathsf{E}_i \mathsf{v}} \cdot \mathsf{A}_i\RoundBrackets*{\mathsf{E}_i \mathsf{v}}.
\]

To analyses $\norm{\mathsf{v}-\mathsf{R}^\intercal_\mathup{c}\mathsf{v}}_{\widetilde{\mathsf{M}}}$, we inevitably require several assumptions for $\widetilde{\mathsf{M}}$. It is known that $\mathsf{A}\lesssim \widetilde{\mathsf{M}}$ via
\[
    \mathsf{I}-\mathsf{A}^{1/2}\widetilde{\mathsf{M}}^{-1}\mathsf{A}^{1/2}=\RoundBrackets*{\mathsf{I}-\mathsf{A}^{1/2}\mathsf{M}^{-1}\mathsf{A}^{1/2}}\RoundBrackets*{\mathsf{I}-\mathsf{A}^{1/2}\mathsf{M}^{-\intercal}\mathsf{A}^{1/2}}.
\]
However, the target here should be bounded $\widetilde{\mathsf{M}}$ by $\mathsf{A}$. A natural thought is to assume $\mathsf{A}\lesssim \widetilde{\mathsf{M}} \lesssim C_\star \mathsf{A}$, which gives exactly an estimate of $\Cond(\widetilde{\mathsf{M}}^{-1}\mathsf{A})$, where $\widetilde{\mathsf{M}}^{-1}\mathsf{A}$ is a preconditioned operator using only the smoother as a preconditioner. Recalling the algebraic representation in \cref{eq:discrete varia}, by using a basic inequality, we can show that
\[
    \begin{aligned}
         & \sum_{i,j}\widetilde{\kappa}_{i+1/2,j}\RoundBrackets*{q_{i+1,j}-q_{i,j}}^2+\widetilde{\kappa}_{i,j+1/2}\RoundBrackets*{q_{i,j+1}-q_{i,j}}^2                     \\
         & \leq 2\sum_{i,j} \widetilde{\kappa}_{i+1/2,j} \RoundBrackets*{q_{i+1,j}^2+q_{i,j}^2}+\widetilde{\kappa}_{i,j+1/2}\RoundBrackets*{q_{i,j+1}^2+q_{i,j}^2}         \\
         & \leq 2\sum_{i,j}\RoundBrackets*{\widetilde{\kappa}_{i-1/2,j}+\widetilde{\kappa}_{i+1/2,j}+\widetilde{\kappa}_{i,j-1/2}+\widetilde{\kappa}_{i,j+1/2}} q_{i,j}^2,
    \end{aligned}
\]
where boundary cases regarding $i$ and $j$ are omitted for brevity. Thanks to harmonic averages, we can see that
\[
    \max\CurlyBrackets*{\widetilde{\kappa}_{i-1/2,j},\widetilde{\kappa}_{i+1/2,j},\widetilde{\kappa}_{i,j-1/2},\widetilde{\kappa}_{i,j+1/2}}\leq 2 \RoundBrackets*{\widetilde{\kappa}_{i,j}^x+\widetilde{\kappa}_{i,j}^y}.
\]
Hence, there exists a generic positive constant $C_\mathup{s}$ such that
\[
    \mathsf{A} \lesssim C_\mathup{s}\sum_{i=1}^{m_\mathup{c}} \mathsf{E}_i^\intercal \mathsf{S}_i \mathsf{E}_i.
\]
One may argue that the matrix $\mathsf{A}$ here is from discretizing a Dirichlet boundary value problem, while the analysis presented is for the original problem \cref{eq:orgional_equation}. Based on the velocity elimination technique, the matrix $\mathsf{A}$ indeed contains an additional term to the original problem. Fortunately, we can show that this term is a positive diagonal matrix and can also be controlled by $\sum_{i=1}^{m_\mathup{c}} \mathsf{E}_i^\intercal \mathsf{S}_i \mathsf{E}_i$ (ref.~\cite{Ye2023}). Then, We can derive that
\[
    \norm{\mathsf{v}-\mathsf{R}^\intercal_\mathup{c}\mathsf{v}_\mathup{c}}_{\widetilde{\mathsf{M}}}^2\leq C_\star \norm{\mathsf{v}-\mathsf{R}^\intercal_\mathup{c}\mathsf{v}_\mathup{c}}_{\mathsf{A}}^2 \leq C_\star C_\mathup{s} \sum_{i=1}^{m_\mathup{c}} \RoundBrackets*{\mathsf{E}_i\mathsf{v}-\mathsf{E}_i\mathsf{R}^\intercal_\mathup{c}\mathsf{v}_\mathup{c}} \cdot \mathsf{S}_i \RoundBrackets*{\mathsf{E}_i\mathsf{v}-\mathsf{E}_i\mathsf{R}^\intercal_\mathup{c}\mathsf{v}_\mathup{c}}.
\]
Recalling that each column of $\mathsf{R}_\mathup{c}^\intercal$ is from solving a generalized eigenvalue problem $\mathsf{A}_i \mathsf{w}=\lambda \mathsf{S}_i \mathsf{w}$ and $\mathsf{E}_i\mathsf{R}^\intercal_\mathup{c}\mathsf{v}_\mathup{c}$ belongs to the local eigenspace, we can utilize the interpolation property to obtain
\[
    \min_{\mathsf{v}_\mathup{c}\in \Real^{n_\mathup{c}}} \sum_{i=1}^{m_\mathup{c}} \RoundBrackets*{\mathsf{E}_i\mathsf{v}-\mathsf{E}_i\mathsf{R}^\intercal_\mathup{c}\mathsf{v}_\mathup{c}} \cdot \mathsf{S}_i \RoundBrackets*{\mathsf{E}_i\mathsf{v}-\mathsf{E}_i\mathsf{R}^\intercal_\mathup{c}\mathsf{v}_\mathup{c}} = \sum_{i=1}^{m_\mathup{c}} \frac{1}{\lambda_{l_\mathup{c}^i}} \RoundBrackets*{\mathsf{E}_i\mathsf{v}} \cdot \mathsf{A}_i \RoundBrackets*{\mathsf{E}_i\mathsf{v}},
\]
where $\lambda_{l_\mathup{c}^i}$ is the $l_\mathup{c}^i$-st smallest\footnote{Using $0$-based numbering here.} eigenvalue. Finally, a lower bound of $\lambda_{\mathup{min}}(\mathsf{P}^{-1}_\SSSText{TG}\mathsf{A})$ can be achieved by controlling the minimal value of $\CurlyBrackets*{\lambda_{l_\mathup{c}^1},\dots,\lambda_{l_\mathup{c}^{m_\mathup{c}}}}$.

The estimate of $\Cond(\mathsf{P}^{-1}_\SSSText{TG}\mathsf{A})$ could be summarized in the following lemma.
\begin{lemma}\label{lem:B_TG}
    Let $C_\star$ be a positive constant such that $\widetilde{\mathsf{M}}\lesssim C_\star\mathsf{A}$, $C_\lambda$ be $\min\CurlyBrackets*{\lambda_{l_\mathup{c}^{i}}}_{i=1}^{m_\mathup{c}}$. Then, the bound
    \[
        \Cond(\mathsf{P}^{-1}_\SSSText{TG}\mathsf{A}) \leq \max\CurlyBrackets*{C_\mathup{s} \frac{C_\star}{C_\lambda},1}
    \]
    holds, where $C_\mathup{s}$ is a generic positive constant.
\end{lemma}

In \cite{Notay2007}, Notay provided two-side bounds of the eigenvalues of $\mathsf{P}_\SSSText{ITG}^{-1}\mathsf{A}$:
\begin{theorem}\label{thm:Notay}
    If  $\mathsf{M}+\mathsf{M}^\intercal-\mathsf{A}$ and $\mathsf{A}_\mathup{c}$ are positive definite matrices. Inequalities
    \[
        \begin{aligned}
            \lambda_{\mathup{max}}(\mathsf{P}_\SSSText{ITG}^{-1}\mathsf{A}) & \leq \lambda_{\mathup{max}}(\mathsf{P}_\SSSText{TG}^{-1}\mathsf{A})\max\CurlyBrackets*{\lambda_{\mathup{max}}(\mathsf{B}_\mathup{c}^{-1}\mathsf{A}_\mathup{c}),1}\ \text{and} \\
            \lambda_{\mathup{min}}(\mathsf{P}_\SSSText{ITG}^{-1}\mathsf{A}) & \geq \lambda_{\mathup{min}}(\mathsf{P}_\SSSText{TG}^{-1}\mathsf{A})\min\CurlyBrackets*{\lambda_{\mathup{min}}(\mathsf{B}_\mathup{c}^{-1}\mathsf{A}_\mathup{c}),1}
        \end{aligned}
    \]
    hold.
\end{theorem}
Investigating recently developed convergence theories for inexact two-grid preconditioners, e.g.~\cite{Xu2022,Xu2022a}, is beyond the scope of the paper. An immediate result of \cref{thm:Notay} is an estimate of the condition number of $\mathsf{P}_\SSSText{ITG}^{-1}\mathsf{A}$ as
\[
    \Cond(\mathsf{P}_\SSSText{ITG}^{-1}\mathsf{A}) \leq \Cond(\mathsf{P}_\SSSText{TG}^{-1}\mathsf{A})\frac{\max\CurlyBrackets*{\lambda_{\mathup{max}}(\mathsf{B}_\mathup{c}^{-1}\mathsf{A}_\mathup{c}),1}}{\min\CurlyBrackets*{\lambda_{\mathup{min}}(\mathsf{B}_\mathup{c}^{-1}\mathsf{A}_\mathup{c}),1}}.
\]

Comparing \cref{eq:B_c} with \cref{eq:B_TG}, we realize that $\mathsf{B}_\mathup{c}$ is an exact two-grid preconditioner for $\mathsf{A}_\mathup{c}$. Hence, two results we mentioned for $\mathsf{P}_\SSSText{TG}^{-1}\mathsf{A}$ can also be applied here, i.e.,
\[
    \lambda_{\mathup{max}}(\mathsf{B}_\mathup{c}^{-1}\mathsf{A}_\mathup{c})=1\ \text{and}\ \frac{1}{\lambda_{\mathup{min}}(\mathsf{B}_\mathup{c}^{-1}\mathsf{A}_\mathup{c})}=\max_{\mathsf{v}_\mathup{c}\in \Real^{n_\mathup{c}}\setminus\CurlyBrackets*{\mathsf{0}}}\min_{\mathsf{v}_\mathup{cc}\in \Real^{n_\mathup{cc}}\setminus\CurlyBrackets*{\mathsf{0}}} \frac{\norm{\mathsf{v}_\mathup{c}-\mathsf{R}_\mathup{cc}^\intercal \mathsf{v}_\mathup{cc}}^2_{\widetilde{\mathsf{M}_\mathup{c}}}}{\norm{\mathsf{v}_\mathup{c}}_{\mathsf{A}_\mathup{c}}^2},
\]
where $\widetilde{\mathsf{M}_\mathup{c}}\coloneqq \mathsf{M}_\mathup{c}^\intercal\RoundBrackets*{\mathsf{M}_\mathup{c}+\mathsf{M}_\mathup{c}^\intercal-\mathsf{A}_\mathup{c}}^{-1}\mathsf{M}_\mathup{c}$. Similarly, for a coarse-coarse element $K_\mathup{cc}^i\in \mathcal{T}_\mathup{cc}$, let $\mathsf{A}_{\mathup{c},i}$ and $\mathsf{S}_{\mathup{c},i}$ be two matrices of the bilinear form of the left-hand and right-hand of \cref{eq:spectral cc} by choosing normalized eigenvectors of \cref{eq:spepb} as bases. We have mentioned that each $\mathsf{S}_{\mathup{c},i}$ is exactly an identity matrix. Meanwhile, we can write an algebraic partition of unity on the coarse grid as $\mathsf{I}=\mathsf{E}_{\mathup{c},1}^\intercal\mathsf{E}_{\mathup{c},1}+\dots+\mathsf{E}_{\mathup{c},m_\mathup{cc}}^\intercal\mathsf{E}_{\mathup{c},m_\mathup{cc}}$, where each $\mathsf{E}_{\mathup{c},i}$ is a binary matrix. Based on the same reason for the relation $\sum_{i=1}^{m_\mathup{c}}\mathsf{E}_i^\intercal \mathsf{A}_i\mathsf{E}_i \lesssim \mathsf{A}$, it can be shown that
\[
    \sum_{i=1}^{m_\mathup{cc}} \mathsf{E}_{\mathup{c},i}^\intercal \mathsf{A}_{\mathup{c},i}\mathsf{E}_{\mathup{c},i} \lesssim \mathsf{A}_\mathup{c}.
\]
Noting that $\widetilde{\mathsf{M}_\mathup{c}}\lesssim \lambda_{\mathup{max}}(\widetilde{\mathsf{M}_\mathup{c}}) \mathsf{I}$, we have
\[
    \begin{aligned}
        \norm{\mathsf{v}_\mathup{c}-\mathsf{R}^\intercal_\mathup{cc}\mathsf{v}_\mathup{cc}}_{\widetilde{\mathsf{M}_\mathup{c}}}^2 & \leq \lambda_{\mathup{max}}(\widetilde{\mathsf{M}_\mathup{c}})\norm{\mathsf{v}_\mathup{c}-\mathsf{R}^\intercal_\mathup{cc}\mathsf{v}_\mathup{cc}}_{\mathsf{I}}^2                                                                                                                                                                                                   \\
                                                                                                                                  & = \lambda_{\mathup{max}}(\widetilde{\mathsf{M}_\mathup{c}}) \sum_{i=1}^{m_\mathup{cc}}\mathsf{E}_{\mathup{c},i}\RoundBrackets*{\mathsf{v}_\mathup{c}-\mathsf{R}^\intercal_\mathup{cc}\mathsf{v}_\mathup{cc}}\cdot \mathsf{S}_{\mathup{c},i}\mathsf{E}_{\mathup{c},i}\RoundBrackets*{\mathsf{v}_\mathup{c}-\mathsf{R}^\intercal_\mathup{cc}\mathsf{v}_\mathup{cc}}.
    \end{aligned}
\]
Therefore, we can cook up an estimate for $\lambda_{\mathup{min}}(\mathsf{B}_\mathup{c}^{-1}\mathsf{A}_\mathup{c})$, which yields the main theorem of this section.
\begin{theorem}\label{thm:main}
    Let $C_\star$, $C_\lambda$ be defined in \cref{lem:B_TG}, the constant $C_{\mathup{c},\lambda}$ be $\min\CurlyBrackets*{\lambda_{l^i_\mathup{cc}}}_{i=1}^{m_\mathup{cc}}$. If $\mathsf{M}+\mathsf{M}^\intercal-\mathsf{A}$, $\mathsf{M}_\mathup{c}+\mathsf{M}^\intercal_\mathup{c}-\mathsf{A}_\mathup{c}$, $\mathsf{A}_\mathup{c}$ and $\mathsf{A}_\mathup{cc}$ are all positive definite matrices, then
    \begin{equation} \label{eq:esti P_ITG}
        \Cond(\mathsf{P}_\SSSText{ITG}^{-1}\mathsf{A}) \leq \max\CurlyBrackets*{C_\mathup{s} \frac{C_\star \lambda_{\mathup{max}}(\widetilde{\mathsf{M}_\mathup{c}})}{C_\lambda C_{\mathup{c},\lambda}},1},
    \end{equation}
    where $C_\mathup{s}$ is a generic positive constant.
\end{theorem}

In \cref{eq:esti P_ITG}, the constants in the numerator are related to the smoothers $\mathsf{M}$ and $\mathsf{M}_\mathup{c}$, while the constants in the denominator are determined by dimensions of local coarse and coarse-coarse spaces. We admit the current form of \cref{thm:main} is far from completeness. For example, the asymptotic rates of $C_\lambda$ and $C_{\mathup{c},\lambda}$ are not depicted. However, our experiments hint that using a modestly large number of eigenvectors is adequate to achieve robustness. Therefore, for \cref{thm:main}, we pursue theoretical guidance rather than rigorous analysis.

\section{Numerical experiments}\label{sec:num}
In this section, numerical experiments are presented to illustrate the performance of the proposed preconditioner.
We consider the unit cube $(0, 1)^3$ as the computational domain $\Omega$ with uniformly structured meshes, i.e.~$N\times N\times N$ with $h^x=h^y=h^z=1/N$ for all experiments, and $n$ the total number of DoF also the fine elements is $N^3$. We implement the proposed spectral three-grid preconditioner in PETSc \cite{Balay2022a} and leverage DMDA---a module provided by PETSc---for inter-process communications and data management. As described in \cref{fig:grid}, each MPI process owns a coarse-coarse element also a cuboid subdomain, the number of MPI processes hence equals $n_\mathup{cc}$ the total number of coarse-coarse elements. Meanwhile, each coarse-coarse element is divided into $\mathtt{sd}$ parts in the $x$-, $y$- and $z$-direction, which defines the coarse grid $\mathcal{T}_\mathup{c}$ and determines the relation $n_\mathup{c}=\mathtt{sd}^3\times n_\mathup{cc}$. To make the notations more self-explanatory, we take $\mathtt{DoF}$ for the number of fine elements (equals $n$), $\mathtt{proc}$ for the number of MPI processes (equals $n_\mathup{cc}$) in the following tables and figures. For the source term $f$ of the model \cref{eq:orgional_equation}, we place $4$ long singular sources at $4$ corner points on the $x\times y$-plane and an opposite sign long singular source at the middle of the $x\times y$-plane, which mimics the well condition in the SPE10 model (ref.~\cite{Christie2001}).

For smoother $\mathsf{M}$ on $\mathsf{A}$, we will apply $\nu$ times block Jacobi iterations, while the block partition is from the coarse grid $\mathcal{T}_\mathup{c}$. Similarly, for $\mathsf{M}_\mathup{c}$, we will take $\nu_\mathup{c}$ times block Jacobi iterations w.r.t.~the coarse-coarse grid $\mathcal{T}_\mathup{cc}$. Therefore, we can fix only the $1$ ghost layer for all cases, and by contrast, we may use different oversampling layers in overlapping domain decomposition methods (cf.~\cite{Ye2023}). In the set-up phase, we will perform factorizations for the matrices $\mathsf{M}$, $\mathsf{M}_\mathup{c}$, and $\mathsf{A}_\mathup{cc}$. Based on our experience, compared to incomplete Cholesky decompositions, full Cholesky decompositions for $\mathsf{M}$ and $\mathsf{M}_\mathup{c}$ could decrease iteration numbers but lead to a significant increase in computation time. Therefore, for $\mathsf{M}$ and $\mathsf{M}_\mathup{c}$, only incomplete Cholesky factorizations are utilized, while the full LU factorization of $\mathsf{A}_\mathup{cc}$ is facilitated by an external package SuperLU\_DIST \cite{Li2003}.

Solving eigenvalue problems is a crucial ingredient in our method, and we employ SLEPc \cite{Hernandez2005}---a companion eigensolver package for PETSc---in the implementation. There are two different strategies in determining dimensions of local coarse and coarse-coarse spaces: (1) use predefined $L_\mathup{c}$ and $L_\mathup{cc}$ eigenvectors for all coarse and coarse-coarse elements; (2) set eigenvalue thresholds $B_\mathup{c}$ and $B_\mathup{cc}$, solve enough candidate eigenvectors and then select eigenvectors whose corresponding eigenvalues are below $B_\mathup{c}$ and $B_\mathup{cc}$ accordingly. The first approach is relatively easy to implement, while the second is based on \cref{thm:main}. The first strategy tends to result in large dimensions of $\mathsf{A}_\mathup{c}$ and $\mathsf{A}_\mathup{cc}$ but balances computing loads among all processes perfectly, while the second is the opposite. We will present a careful comparison of those two strategies later.

We choose the PETSc GMRES method with default parameters to solve the fine system \cref{eq:fine_system} with a different relative tolerance $10^{-6}$, and iteration counts are indicated by $\mathtt{iter}$. The HPC cluster in which we run the program is interconnected with an InfiniBand network, and each node is equipped with dual AMD\textsuperscript{\textregistered}~7452 CPUs ($64$ CPU cores in total) and 256 GB of memory. We will frequently write the elapsed wall time $\mathtt{time}$ and iteration count $\mathtt{iter}$ of an experiment in a format of $\mathtt{time}(\mathtt{iter})$, where one effective digit is reserved in the unit of seconds for $\mathtt{time}$. Our codes are hosted on GitHub\footnote{\href{https://github.com/Laphet/PreMixFEM.git}{https://github.com/Laphet/PreMixFEM.git}}.

\subsection{Scalability tests}
\begin{figure}[!ht]
    \centering
    \begin{subfigure}[b]{0.49\textwidth}
        \includegraphics[width=\textwidth]{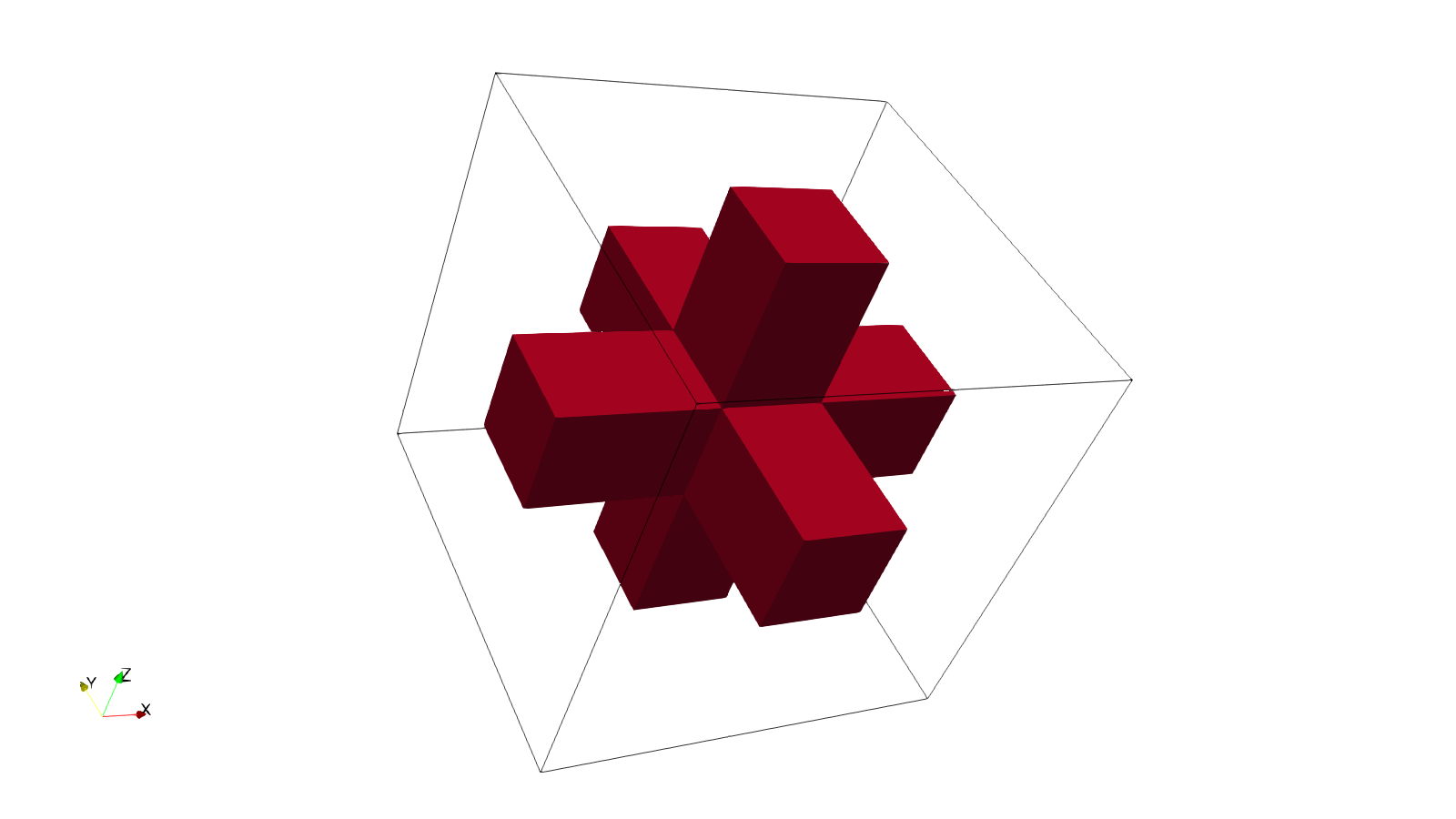}
        \caption{}\label{fig:scal a}
    \end{subfigure}
    \begin{subfigure}[b]{0.49\textwidth}
        \includegraphics[width=\textwidth]{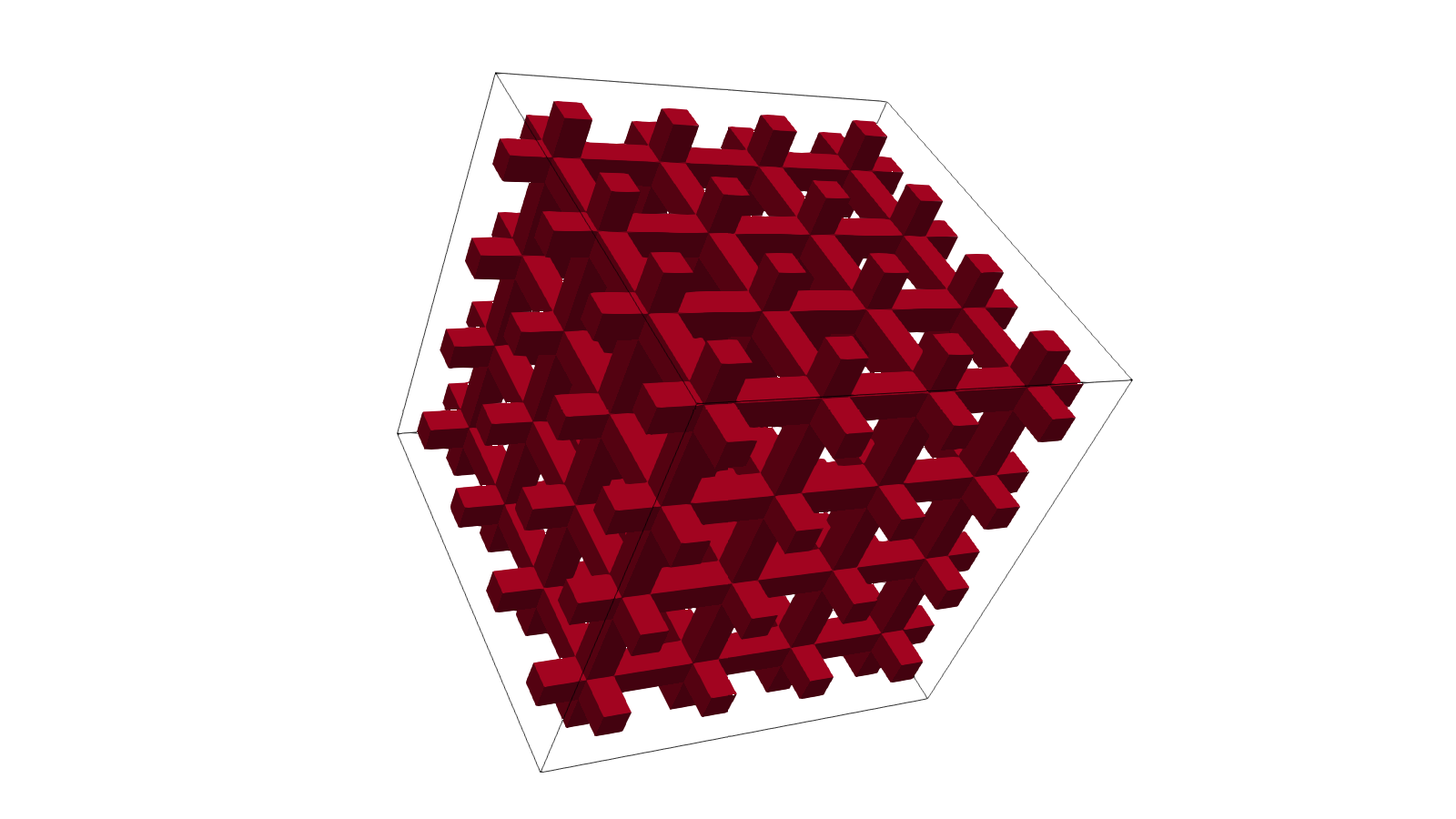}
        \caption{}\label{fig:scal b}
    \end{subfigure}
    \begin{subfigure}[b]{\textwidth}
        \includegraphics[width=\textwidth]{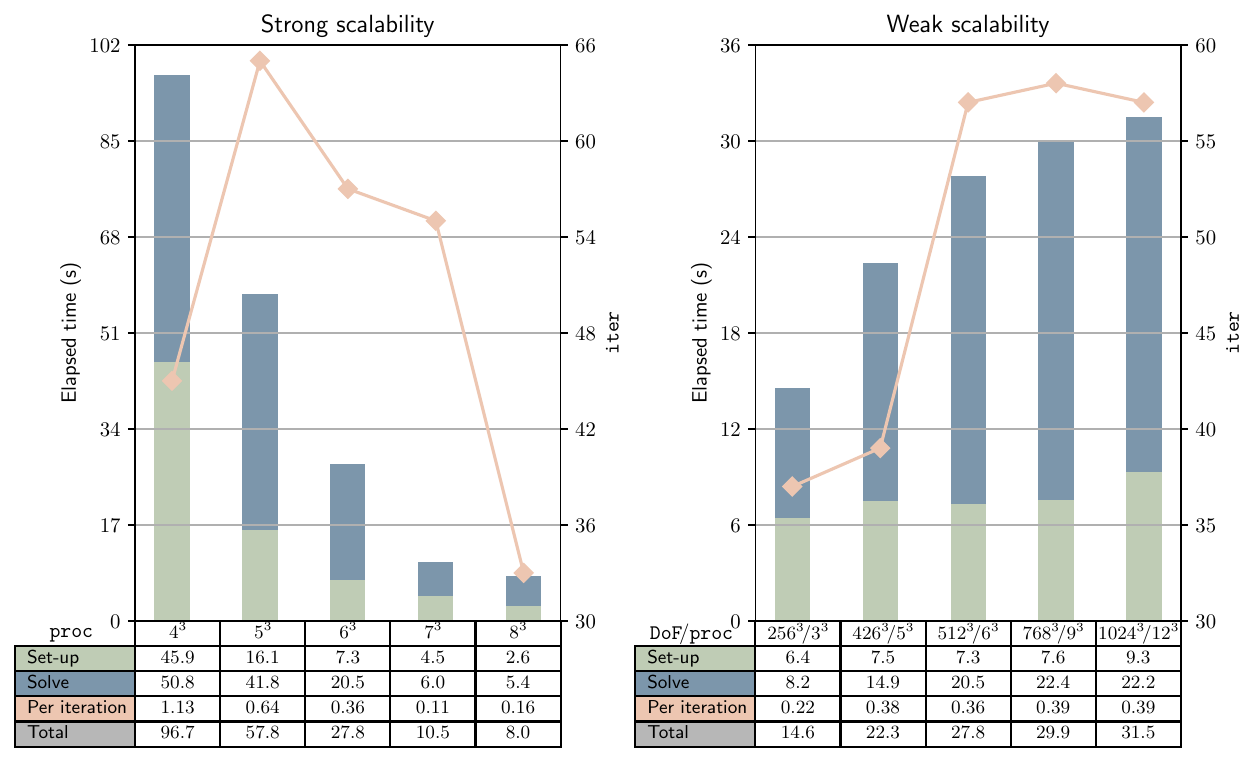}
        \caption{}\label{fig:scal c}
    \end{subfigure}

    \caption{(a) The periodic cell configuration used for the scalability tests. (b) An illustration of a domain that consists of $4\times 4\times 4$ cells, where long channels across the domain. (c)-\textbf{left} the results of strong scalability tests, where the DoF is fixed as $512^3$ and the number of MPI processes varies from $\CurlyBrackets{4^3, 5^3, 6^3, 7^3, 8^3}$. (c)-\textbf{right} the results of weak scalability tests, where the ratio of $\mathtt{DoF}$ and $\mathtt{proc}$ is fixed around $85^3$.}
\end{figure}

In this subsection, we fix $(\mathtt{sd},\nu,\nu_\mathup{c})=(7, 1, 1)$ and adopt the first strategy to construct coarse and coarse-coarse spaces with $(L_\mathup{c}, L_\mathup{cc})=(4, 4)$. A domain $\Omega$ consists of periodically duplicated cells, and a cell consists of $8\times 8 \times 8$ fine elements. Hence, a domain contains $32\times 32 \times 32$ cells if $\mathtt{DoF}=256^3$. We set an isotropic $\mathbb{K}$ as $\mathbb{K}(\bm{x})=\kappa(\bm{x})\mathsf{I}$. The coefficient $\kappa$ takes values from $\CurlyBrackets{1, 10^6}$, and the region of $\kappa(\bm{x})=10^6$ inside a cell is visualized in \cref{fig:scal a}. We demonstrate a domain constructed from $4\times 4 \times 4$ cells in \cref{fig:scal b}, and note that all channels are placed with a high coefficient ($10^6$) and connected, which is considered as a challenging case for traditional solvers (see \cite{Dryja1996}). Note that from the setting, the dimension of the square matrix $\mathsf{A}_\mathup{cc}$ is $\mathtt{proc}\times L_\mathup{cc}$, while for $\mathsf{A}_\mathup{c}$ it is $\mathtt{proc}\times \mathtt{sd}^3\times L_\mathup{c}$, which is $343$ times larger. In our previous work \cite{Ye2023}, we directly factorized $\mathsf{A}_\mathup{c}$, which accounts for an unsatisfactory scalability performance. We shall witness that the three-grid setting improves scalability significantly.

For strong scalability, we set $\mathtt{DoF}=512^3$ and choose $\mathtt{proc}$ from $\CurlyBrackets{4^3, 5^3, 6^3, 7^3, 8^3}$. The results are shown in the left part of \cref{fig:scal c}, which contains a graph and a table. A solution process is divided into ``Set-up'' and ``Solve'' phases that are colored differently. Recall that ``Set-up'' contains solving eigenvalue problems and factorizing matrices. Because we only involve a direct solver in the LU decomposition of $\mathsf{A}_\mathup{cc}$, as expected, the time consumed on ``Set-up'' for $\mathtt{proc}=8^3$ is around $6\%$ of $\mathtt{proc}=4^3$, which is a satisfactory decrease. As a comparison, in our previous work \cite{Ye2023}, the time consumed on the factorization on the coarse level can increase. We also note that $\mathtt{iter}$ does not strictly decrease w.r.t.~$\mathtt{proc}$ in contrast to what we have observed in \cite{Ye2023}. The reason could be that the coarse grid is not aligned with cells, while it is true in \cite{Ye2023}. To exclude the influence of different iteration numbers, we calculate the time per iteration in the table. In the ``Per iteration'' row, we can see that the time decreases by a factor of $10$ for $\mathtt{proc}$ going from $4^3$ to $7^3$ but slightly increases for $\mathtt{proc}=8^3$, and this increase here may be attributed internode communications.

For weak scalability, we choose $\mathtt{DoF}$ from $\CurlyBrackets{256^3, 426^3, 512^3, 768^3, 1024^3}$ and $\mathtt{proc}$ from $\CurlyBrackets{3^3, 5^3, 6^3, 9^3, 12^3}$ accordingly, while the ratio of $\mathtt{DoF}$ and $\mathtt{proc}$ is fixed around $85^3$. The results are shown in the right part of \cref{fig:scal c}. We can see the time consumed on ``Set-up'' is stable w.r.t.~$\mathtt{proc}$ from $4^3$ to $9^3$, and a deterioration of $(9.3/7.6-1.0\approx 12\%)$ is observed for $\mathtt{proc}=12^3$. From the ``Per iteration'' row, the time records of $\mathtt{proc}=5^3$ to $12^3$ are quite close. We may postulate that the $15\%$ deterioration is due to the direct solver. We could also claim that the current three-grid framework is capable of handling problems of $1$ billion DoF or below. However, to extend the proposed method to a larger problem size, truly multigrid and sophisticated load-balancing techniques are crucial. Because there are no inter-node communications needed for $\mathtt{proc}=3^3$, one may notice that the time, in this case, is considerably less than others.

\subsection{Robustness tests on a fractured medium}\label{subsec:robustness}
\begin{figure}[!ht]
    \centering
    \includegraphics[width=3.5in]{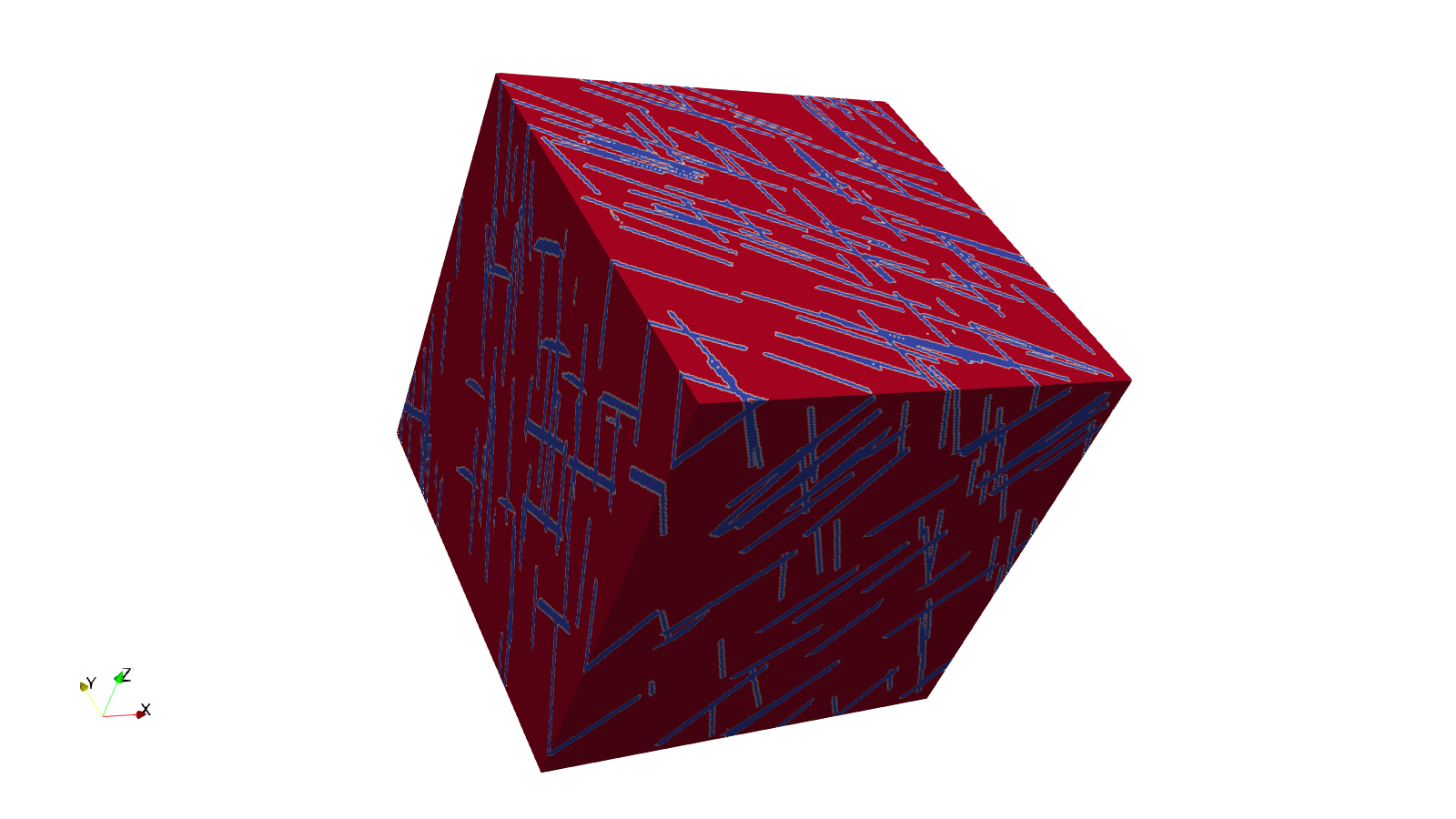}
    \caption{A fractured medium model, where red and blue regions are referred to as matrix (place $1$ as the conductivity) and fractures (place $10^\mathtt{cr}$ as the conductivity) respectively.}
    \label{fig:fractured}
\end{figure}

Modeling flow through fractured porous media is a challenging but important task in reservoir simulation today. Most works treat fractures as planes embedded in the matrix, and the discretization of the domain should be carefully handled (ref.~\cite{Berre2018}). We here simply discretize the medium with a fine uniform mesh and assign large permeability for fine elements belonging to fractures. In \cref{fig:fractured}, we demonstrate a fractured medium, where the red and blue regions correspond to the matrix and fractures, respectively. We again set an isotropic $\mathbb{K}$ as $\mathbb{K}(\bm{x})=\kappa(\bm{x})\mathsf{I}$, where $\kappa=1.0$ for matrix and $\kappa=10^\mathtt{cr}$ for fractures. The original resolution of the medium in \cref{fig:fractured} is $256\times 256 \times 256$, and we periodically duplicate it in the domain $\Omega$ so that $\mathtt{DoF}=512^3$.

In this subsection, we adopt the first strategy to construct coarse and coarse-coarse spaces and fix $(\mathtt{proc}, \mathtt{sd}, \nu, \nu_\mathup{c})=(6^3,7, 1, 1)$. We focus on examining the robustness of the proposed preconditioner w.r.t.~contrast ratios and hence choose $\mathtt{cr}$ from $\CurlyBrackets{1,2,3,4,5,6}$. As a comparison, performances of two off-the-sheaf preconditioners ASM (Additive Schwarz Method) and GAMG (General Algebraic Multigrid) from PETSc are also recorded. Note that we did not tailor special parameters or options for those preconditioners and just used them in a black-box manner. The results are presented in \cref{tab:fracture1}, where we fix $L_\mathup{c}=4$ and choose $L_\mathup{cc}$ from $\CurlyBrackets{4, 8, 11, 17}$ for the proposed preconditioner. Here the values $\CurlyBrackets{4, 8, 11, 17}$ are heuristically determined by the ladder of the eigenvalue sequence
\[
    \CurlyBrackets*{0,\ 1,\ 1,\ 1,\quad 2,\ 2,\ 2,\ 3,\quad 4,\ 4,\ 4,\quad 5,\ 5,\ 5,\ 5,\ 5,\ 5,\quad 6,\dots}
\]
of the Laplace operator.

\begin{table}[htbp]
    \footnotesize
    \caption{For the fractured medium model, records of elapsed wall time and iteration numbers of ASM, GAMG, and the proposed preconditioner w.r.t.~different contrast ratios, where $(*)$ means the solver does not converge in $1000$ iterations and $(k^\Diamond)$ means there is a BREAKDOWN error thrown at the step $k$.}\label{tab:fracture1}
    \begin{center}
        \begin{tabular}{|c|c|c|c|c|c|c|} \hline
            Preconditioner                    & $\mathtt{cr}=1$ & $\mathtt{cr}=2$ & $\mathtt{cr}=3$ & $\mathtt{cr}=4$ & $\mathtt{cr}=5$     & $\mathtt{cr}=6$     \\ \hline
            ASM                               & $135.8(*)$      & $134.5(*)$      & $136.8(*)$      & $135.8(*)$      & $137.1(*)$          & $135.1(*)$          \\ \hline
            GAMG                              & $15.9(14)$      & $19.7(23)$      & $26.8(43)$      & $54.8(117)$     & $45.5(90^\Diamond)$ & $46.1(90^\Diamond)$ \\ \hline
            $L_\mathup{c}=4,L_\mathup{cc}=4$  & $21.9(38)$      & $32.8(70)$      & $40.6(93)$      & $53.9(132)$     & $54.8(137)$         & $55.4(139)$         \\ \hline
            $L_\mathup{c}=4,L_\mathup{cc}=8$  & $21.0(31)$      & $26.9(49)$      & $35.1(72)$      & $41.9(89)$      & $49.1(110)$         & $49.7(112)$         \\ \hline
            $L_\mathup{c}=4,L_\mathup{cc}=11$ & $20.5(28)$      & $26.6(44)$      & $31.9(60)$      & $38.2(75)$      & $43.2(88)$          & $43.4(89)$          \\ \hline
            $L_\mathup{c}=4,L_\mathup{cc}=17$ & $22.6(26)$      & $28.0(39)$      & $31.7(49)$      & $35.3(56)$      & $37.6(62)$          & $37.3(63)$          \\ \hline
        \end{tabular}
    \end{center}
\end{table}

From \cref{tab:fracture1}, we can see that ASM performs poorly for this strongly heterogeneous model, and the solver cannot converge in $1000$ iterations even for $\mathtt{cr}=1$. As $\mathtt{cr}$ increases, the performance of GAMG gradually degrades, i.e., the computation time and iteration numbers all increase. Moreover, we observed BREAKDOWN errors for $\mathtt{cr}=5$ and $6$. For the proposed preconditioner, we can notice that it is quite stable w.r.t.~contrast ratios. The iteration numbers cannot be completely irrelevant with $\mathtt{cr}$ but are well-controlled. Importantly, we did not experience any BREAKDOWN errors in all test cases of the proposed preconditioner. We can also see that including more eigenvectors in the coarse-coarse space greatly stabilizes the performance. For example, the computation time from $\mathtt{cr}=1$ to $\mathtt{cr}=6$ increases only by a factor of $68\%$ for $L_\mathup{cc}=17$, while this number is $152\%$ for $L_\mathup{cc}=4$.

We also carried out several experiments to investigate the influence of $L_\mathup{c}$ on the dimension of local coarse spaces, and the results are shown in \cref{tab:fracture2}, where $\mathtt{cr}=6$ is set. We can see that a larger $L_\mathup{c}$ could improve the performance, but not significantly compared to $L_\mathup{cc}$. Specifically, if $L_\mathup{cc}=4$, the computation time increases by a factor of $32\%$ from $L_\mathup{c}=4$ to $17$; if $L_\mathup{c}=4$, the time decreases by a factor of $32\%$ from $L_\mathup{cc}=4$ to $17$. Although the minimal iteration number occurs at $(L_\mathup{c}, L_\mathup{cc})=(17, 17)$, the shortest time is found at $(L_\mathup{c}, L_\mathup{cc})=(4, 17)$. The combination effect of $L_\mathup{c}$ and $L_\mathup{cc}$ is complicated to be disentangled, while it is advised to enlarge $L_\mathup{cc}$ rather than $L_\mathup{c}$ first for improving the final performance.

\begin{table}[htbp]
    \footnotesize
    \caption{For the fractured medium model, records of elapsed wall time and iteration numbers w.r.t.~different $L_\mathup{c}$ and $L_\mathup{cc}$.}\label{tab:fracture2}
    \begin{center}
        \begin{tabular}{|c|c|c|c|c|}\hline
            $L_\mathup{cc}$ & $L_\mathup{c}=4$ & $L_\mathup{c}=8$ & $L_\mathup{c}=11$ & $L_\mathup{c}=17$ \\ \hline
            $4$             & $55.4(139)$      & $60.0(133)$      & $63.1(129)$       & $72.8(126)$       \\ \hline
            $8$             & $49.7(112)$      & $55.3(114)$      & $57.3(110)$       & $67.2(106)$       \\ \hline
            $11$            & $44.3(89)$       & $48.3(92)$       & $51.6(90)$        & $60.5(87)$        \\ \hline
            $17$            & $37.3(63)$       & $40.6(61)$       & $42.5(59)$        & $50.6(58)$        \\ \hline
        \end{tabular}
    \end{center}
\end{table}

\subsection{Tests on the second strategy}
It is impossible for SLEPc the eigensolver that we used, to set a target priorly and then find all eigenvectors whose corresponding eigenvalues are below this target. We hence first take a relatively large number into the solver routine, which will produce this number of eigenvectors whose eigenvalues are in ascending order. We can then collect a subset of those eigenvectors and utilize them to construct the local coarse or coarse-coarse space. According to results in \cref{subsec:robustness}, dimensions of the coarse-coarse space contribute more to the final performance compared with the coarse space. Therefore, we only take the second strategy to construct local coarse-coarse spaces and will fix $L_\mathup{c}=4$ in the following experiments. On each coarse-coarse element, we compute $20$ candidate eigenvectors and select those whose corresponding eigenvalues are below $B_\mathup{cc}$. The construction of the coefficient field follows the same manner as in \cref{subsec:robustness}. For other parameters, we set $(\mathtt{proc}, \mathtt{sd}, \nu, \nu_\mathup{c})=(6^3,7, 1, 1)$. The results are presented in \cref{tab:second approach}, where the dimensions of the coarse-coarse spaces are also tabulated. In the last row of \cref{tab:second approach}, we set $B_\mathup{cc}=10^{12}$ as a reference, which is basically back to the first strategy with $L_\mathup{cc}=20$.

\begin{table}[htbp]
    \footnotesize
    \caption{Records of elapsed wall time, iteration numbers and dimensions of coarse-coarse spaces w.r.t.~different $\mathtt{cr}$ and $B_\mathup{cc}$.} \label{tab:second approach}
    \begin{center}
        \begin{tabular}{|c|c|c|c|c|c|c|}\hline
            $B_\mathup{cc}$                      & $\mathtt{cr}=1$ & $\mathtt{cr}=2$ & $\mathtt{cr}=3$ & $\mathtt{cr}=4$ & $\mathtt{cr}=5$ & $\mathtt{cr}=6$ \\ \hline
            \multirow{2}{*}{$2.0\times 10^{-4}$} & $23.3(43)$      & $33.4(71)$      & $37.3(73)$      & $40.4(72)$      & $40.4(74)$      & $37.9(75)$      \\ \cline{2-7}
                                                 & $849$           & $895$           & $1748$          & $2125$          & $2381$          & $2381$          \\ \hline
            \multirow{2}{*}{$4.0\times 10^{-4}$} & $21.1(32)$      & $27.3(47)$      & $32.9(53)$      & $33.5(58)$      & $37.4(62)$      & $38.0(63)$      \\ \cline{2-7}
                                                 & $1491$          & $1723$          & $2742$          & $3169$          & $3222$          & $3222$          \\ \hline
            \multirow{2}{*}{$8.0\times 10^{-4}$} & $18.9(27)$      & $29.1(39)$      & $30.7(46)$      & $35.5(53)$      & $38.8(60)$      & $39.8(60)$      \\ \cline{2-7}
                                                 & $2679$          & $3438$          & $4109$          & $4218$          & $4220$          & $4221$          \\ \hline
            \multirow{2}{*}{$10^{12}$}           & $23.3(24)$      & $29.3(37)$      & $32.6(46)$      & $37.1(53)$      & $39.4(60)$      & $40.4(60)$      \\ \cline{2-7}
                                                 & $4320$          & $4320$          & $4320$          & $4320$          & $4320$          & $4320$          \\ \hline
        \end{tabular}
    \end{center}
\end{table}

We can notice that for $B_\mathup{c}\in \CurlyBrackets{2.0, 4.0, 8.0}\times 10^{-4}$, the dimension of the resulting coarse-coarse space increases with $\mathtt{cr}$. It reveals a pattern that eigenvalues of the spectral problem concentrate toward $0$ as increasing the contrast ratio, which may explain the necessity of incorporating enough eigenvectors for high contrast problems. But fortunately, the tendency gets saturated as the contrast ratio becomes large, i.e., variances of dimensions for $\mathtt{cr} \in \CurlyBrackets{4, 5, 6}$ are noticeably smaller compared with $\mathtt{cr}\in\CurlyBrackets{1, 2, 3}$. This observation also explains that the second strategy is more effective for small contrast problems. We can find the following fact to support it: the computation time of $B_\mathup{cc}=8.0\times 10^{-4}$ at $\mathtt{cr}=1$ is equal to around $81\%$ of $B_\mathup{cc}=10^{12}$, while this value is $98\%$ at $\mathtt{cr}=6$. Moreover, we can also conclude that the second strategy exhibits less advantage or even disadvantage over the first strategy for high contrast cases, e.g., the best computation time record in \cref{tab:fracture2} is slightly less than \cref{tab:second approach} at $\mathtt{cr}=6$. To fully leverage the DoF reduction from the second strategy, a load re-balancing procedure should be emphasized and exploited, while it is missed in our current implementation.

\subsection{Comparison with the exact preconditioner}
We have emphasized repeatedly the defect of the exact two-grid preconditioner, and in this subsection, we will provide several concrete experiments to support the point. We still use the same way in \cref{subsec:robustness} to produce the permeability field and adopt the first strategy to building local coarse and coarse-coarse spaces. The results are shown in \cref{tab:exact}, where other parameters are set as $(\mathtt{proc}, L_\mathup{c}, L_\mathup{cc}, \nu, \nu_\mathup{c})=(6^3, 4, 8, 1, 1)$. In the last column of \cref{tab:exact}, we also calculate the dimensions of $W_\mathup{c}$ and $W_\mathup{cc}$ via $\mathtt{proc}\times \mathtt{sd}^3 \times L_\mathup{c}$ and $\mathtt{proc}\times L_\mathup{cc}$.

\begin{table}[htbp]
    \footnotesize
    \caption{Records of elapsed wall time, iteration numbers and dimensions of $W_\mathup{c}$ or $W_\mathup{cc}$ of the exact and inexact preconditioner w.r.t.~different $\mathtt{cr}$ and $\mathtt{sd}$.} \label{tab:exact}
    \begin{center}
        \begin{tabular}{|c|c|c|c|c|c|c|c|c|}\hline
            $\mathtt{sd}$        & Type    & $\mathtt{cr}=1$ & $\mathtt{cr}=2$ & $\mathtt{cr}=3$ & $\mathtt{cr}=4$ & $\mathtt{cr}=5$ & $\mathtt{cr}=6$ & Dim      \\ \hline
            \multirow{2}{*}{$6$} & exact   & $32.0(24)$      & $34.6(33)$      & $37.2(42)$      & $39.4(49)$      & $41.2(50)$      & $41.3(50)$      & $186624$ \\ \cline{2-9}
                                 & inexact & $24.1(32)$      & $31.5(51)$      & $41.7(77)$      & $46.8(91)$      & $53.5(108)$     & $52.3(108)$     & $1728$   \\ \hline
            \multirow{2}{*}{$7$} & exact   & $36.9(22)$      & $39.6(29)$      & $45.9(34)$      & $45.4(37)$      & $47.6(37)$      & $48.8(38)$      & $296352$ \\ \cline{2-9}
                                 & inexact & $21.0(31)$      & $26.9(49)$      & $35.1(72)$      & $41.9(89)$      & $49.1(110)$     & $49.7(112)$     & $1728$   \\ \hline
            \multirow{2}{*}{$8$} & exact   & $45.8(20)$      & $48.0(26)$      & $50.4(30)$      & $50.3(32)$      & $57.0(32)$      & $48.5(32)$      & $442368$ \\ \cline{2-9}
                                 & inexact & $19.7(30)$      & $26.1(46)$      & $37.5(72)$      & $43.6(88)$      & $55.0(110)$     & $53.0(118)$     & $1728$   \\ \hline
        \end{tabular}
    \end{center}
\end{table}

From \cref{tab:exact}, we can notice that iteration numbers decrease as $\mathtt{sd}$ increases in all cases. This phenomenon may be due to the domain decomposition in the implementation that a larger $\mathtt{sd}$ means a smaller coarse element, which complies with the convergence theories of traditional geometric multigrid preconditioners. We can also see that the exact preconditioner is much more robust for contrast ratios than the inexact one. For example, when $\mathtt{sd}=8$ for the exact preconditioner, the iteration number increases by a factor of $60\%$ from $\mathtt{cr}=1$ to $\mathtt{cr}=6$, while this value is $293\%$ for the inexact preconditioner. However, a less iteration number does not always mean better performance. To see this, the records of computation time for the exact preconditioner for $\mathtt{cr}=6$ are longer than the best one in \cref{tab:fracture2}.

\section{Application on two-phase flow problem}\label{sec:spe10}
In this section, we will apply the proposed preconditioner to solve immiscible and incompressible two-phase (water and oil) flow problems. Physical quantities corresponding to different phases are subscripted as $\mathup{o}$ (oil) and $\mathup{w}$ (water). We assume that the fluid is incompressible, and the governing transportation equation for each phase $\alpha \in \CurlyBrackets{\mathup{w}, \mathup{o}}$ is stated as follows:
\begin{equation} \label{eq:mass}
    \phi\frac{\partial S_\alpha}{\partial t} + \Div {\bm{v}}_\alpha = \widetilde{q}_\alpha,
\end{equation}
where $\phi$ is the porosity of the reservoir rock, and for $\alpha \in \CurlyBrackets{\mathup{w}, \mathup{o}}$, $S_\alpha$ is the saturation, $\bm{v}_\alpha$ is the superficial Darcy velocity that is given by Darcy's law, $\widetilde{q}_\alpha$ is the external sources and sinks but divided by phase's density. Darcy's law for two-phase flows is expressed in the following form:
\begin{equation} \label{eq:darcy}
    \bm{v}_\alpha=-\frac{k_\alpha^\mathup{r}}{\mu_\alpha} \mathbb{K} \nabla {p_\alpha}.
\end{equation}
where $k_\alpha^\mathup{r}$ is the relative permeability that depends on the saturation, $\mu_\alpha$ is the viscosity, and $\mathbb{K}$ is the orthotropic permeability of the reservoir rock, and $p_\alpha$ is the pressure. Note that gravity is neglected in \cref{eq:darcy} for brevity. The fact that the two fluids jointly fill the voids implies the relation
\begin{equation}\label{eq:unit}
    S_\mathup{w}+S_\mathup{o}=1.
\end{equation}
For simplicity, we also neglect the capillary pressure and hence have
\begin{equation}\label{eq:pressure diff}
    p_\mathup{o}-p_\mathup{w}=0.
\end{equation}
The system of \cref{eq:mass,eq:darcy,eq:unit,eq:pressure diff} is overdetermined, and we usually assign $p_\mathup{o}$ and $S_\mathup{w}$ as the two primary variables to cancel out other equations. By defining the total velocity
\[
    \bm{v}=\bm{v}_\mathup{w}+\bm{v}_\mathup{o}
\]
and redefining $p\coloneqq p_\mathup{o}$, $S\coloneqq S_\mathup{w}$, \cref{eq:mass,eq:darcy,eq:unit,eq:pressure diff} are reduced into
\begin{align}
    \bm{v}                                                                        & =-\lambda(S)\mathbb{K}\nabla p, \label{eq:v}                                                \\
    \Div \bm{v}                                                                   & =\widetilde{q}\coloneqq \widetilde{q}_\mathup{w}+\widetilde{q}_\mathup{o}, \label{eq:div_v} \\
    \phi\frac{\partial S}{\partial t}+\Div \RoundBrackets*{f_\mathup{w}(S)\bm{v}} & =\widetilde{q}_\mathup{w} \label{eq:trans_S},
\end{align}
where $\lambda\coloneqq k_\mathup{w}^\mathup{r}/\mu_\mathup{w}+k_\mathup{o}^\mathup{r}/\mu_\mathup{o}$ is called the total mobility, and $f_\mathup{w}\coloneqq\RoundBrackets*{k_\mathup{w}^\mathup{r}/\mu_\mathup{w}}/\lambda$ is called the fractional flow function.

\begin{figure}[!ht]
    \centering
    \begin{subfigure}[b]{0.49\textwidth}
        \centering
        \includegraphics[width=\textwidth]{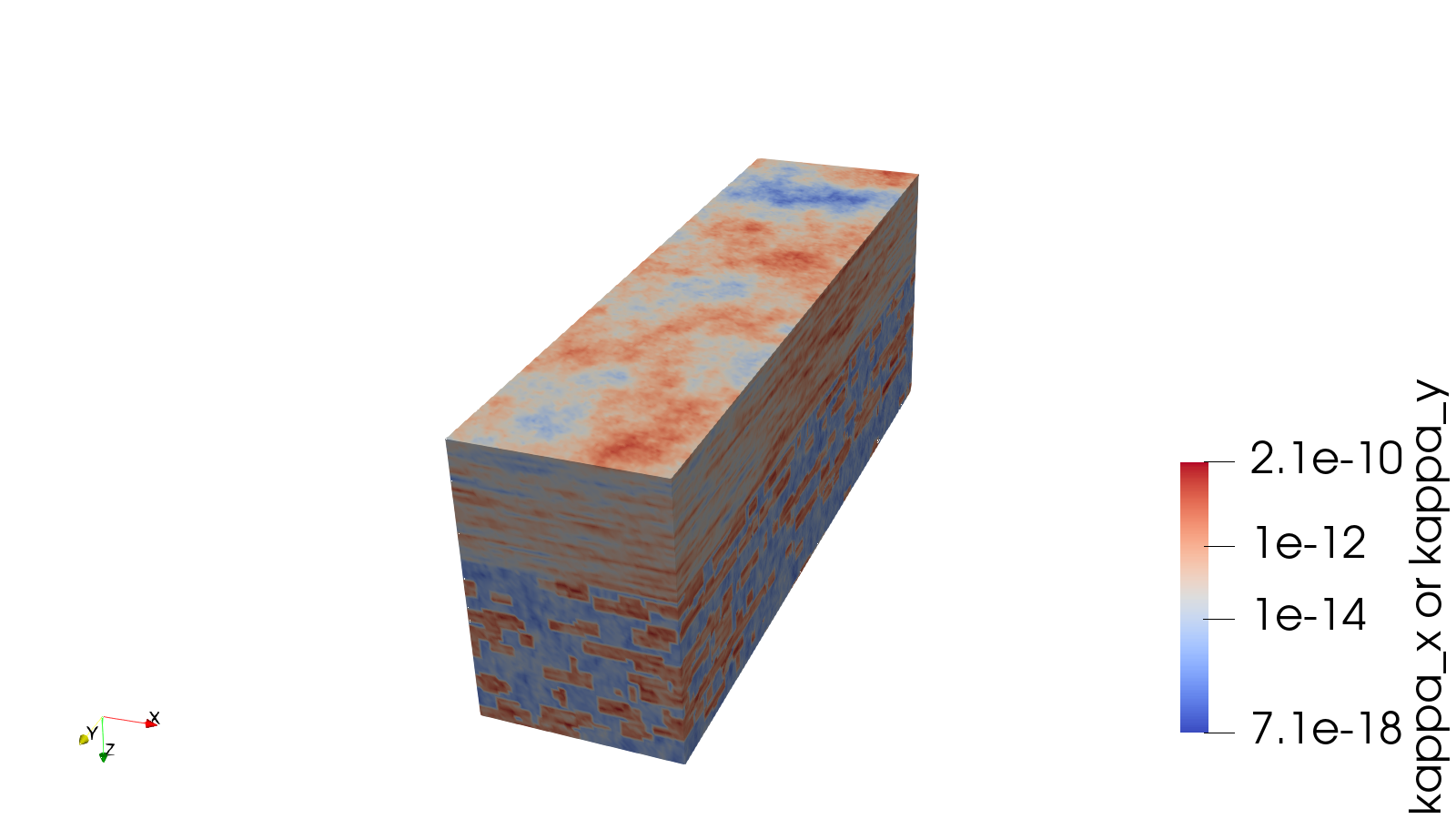}
        \caption{}\label{fig:spe10_perm_x}
    \end{subfigure}
    \begin{subfigure}[b]{0.49\textwidth}
        \centering
        \includegraphics[width=\textwidth]{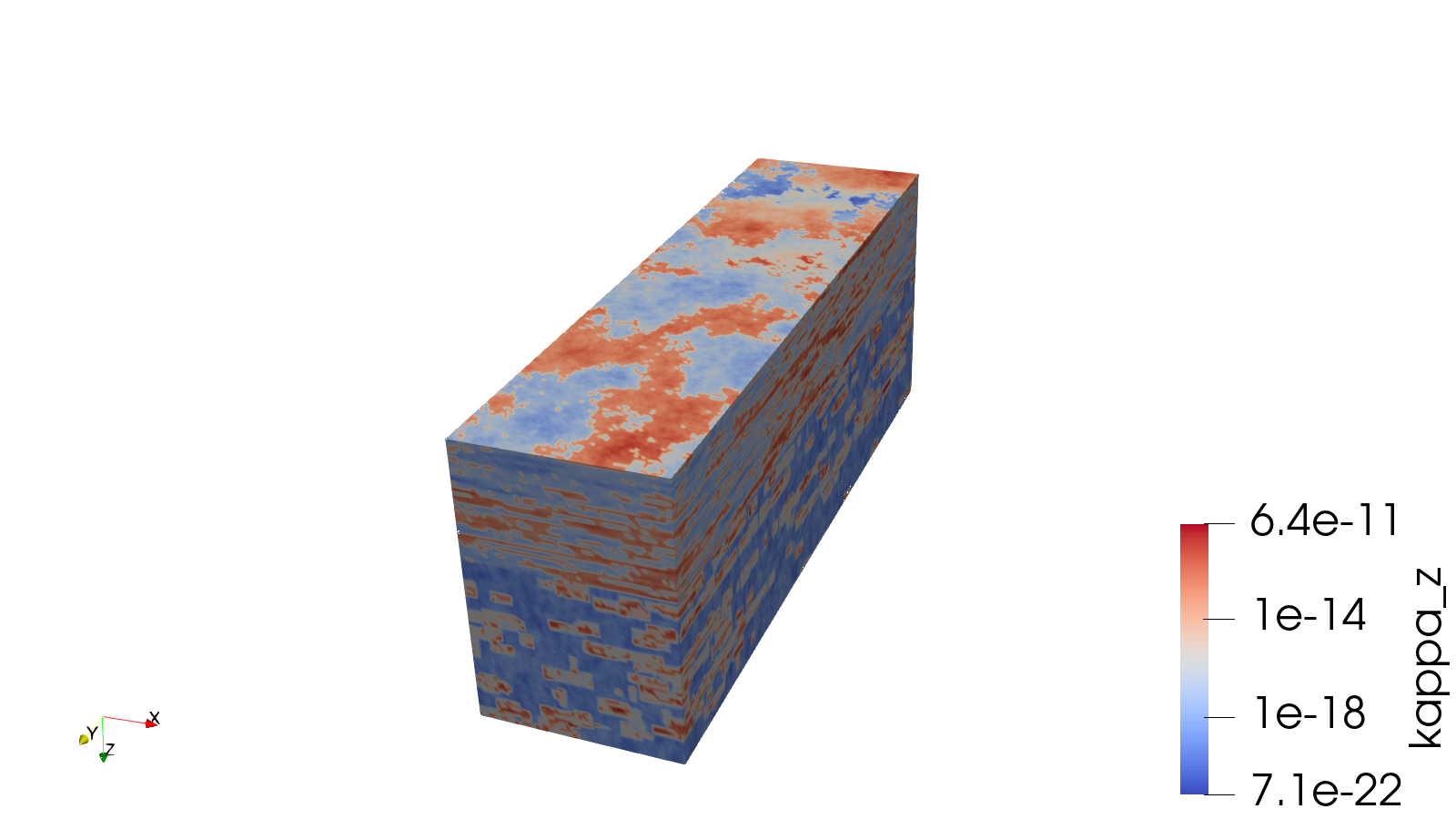}
        \caption{}\label{fig:spe10_perm_z}
    \end{subfigure}
    \begin{subfigure}[b]{0.49\textwidth}
        \centering
        \includegraphics[width=\textwidth]{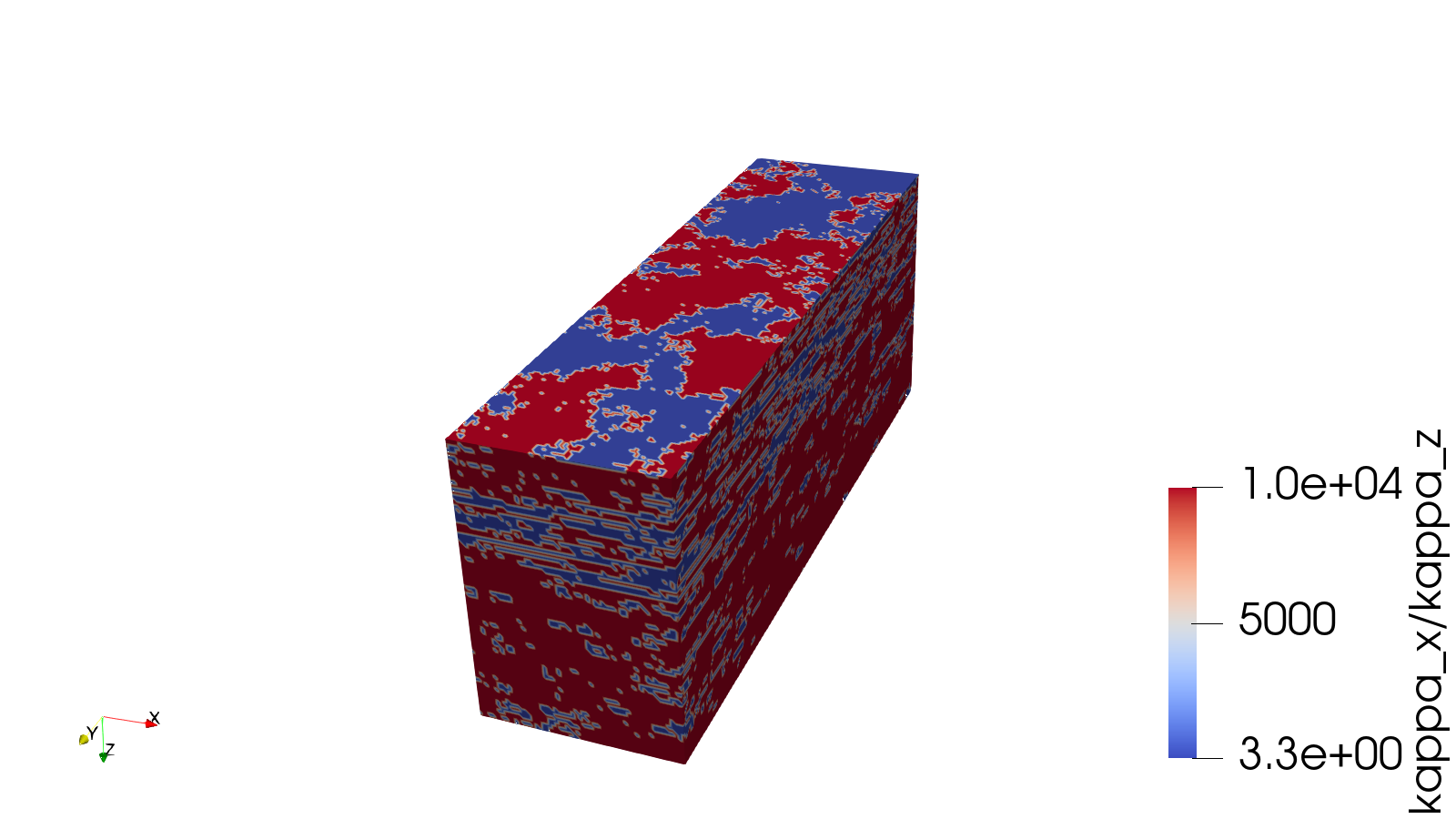}
        \caption{}\label{fig:spe10_ratio}
    \end{subfigure}
    \begin{subfigure}[b]{0.49\textwidth}
        \centering
        \includegraphics[width=\textwidth]{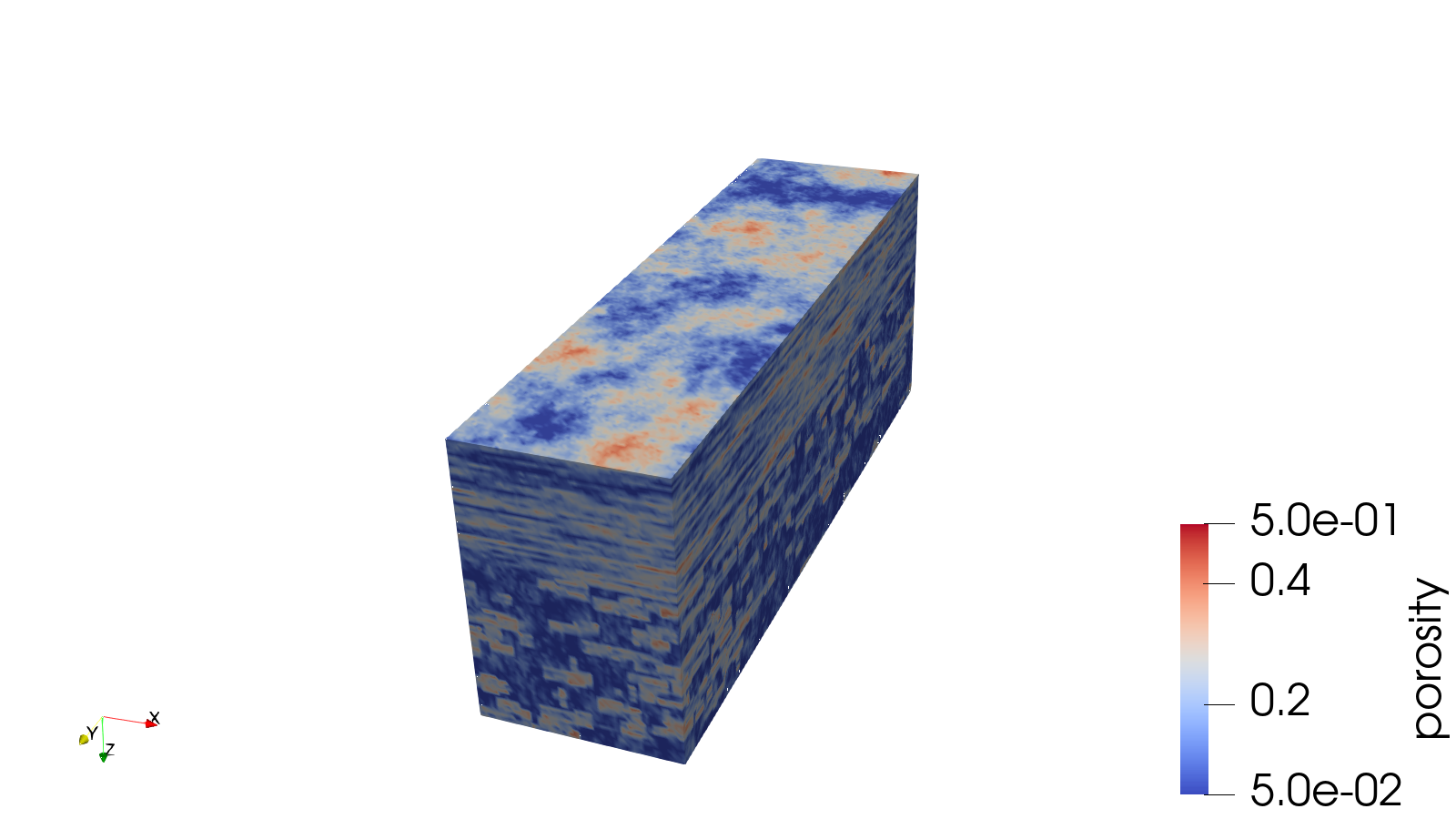}
        \caption{}\label{fig:spe10_poro}
    \end{subfigure}
    \caption{Rock properties of the SPE10 model: (a) $\kappa^x$ or $\kappa^y$, (b) $\kappa^z$, (c) $\kappa^x\!/\kappa^z$, and (d) modified $\phi$.}
\end{figure}

The system of \cref{eq:v,eq:div_v,eq:trans_S} is incomplete; for example, the law of $k^\mathup{r}_\alpha$ w.r.t.~$S_\alpha$ is missing. We introduce the SPE10 model (ref.~\cite{Christie2001}) to achieve a thorough description of the problem. The computational domain of the SPE10 model has a physical dimension of $1200\,\mathup{ft}\times 2200\,\mathup{ft}\times 170\,\mathup{ft}$ and contains $60 \times 220 \times 85$ cells, which yields the dimension of fine elements as $(h^x, h^y, h^z)=(20,10,2)\,\mathup{ft}$. For the permeability field $\mathbb{K}=\Diag(\kappa^x,\kappa^y,\kappa^z)$, $\kappa^x$ is the same as $\kappa^y$ and visualized in \cref{fig:spe10_perm_x}, while $\kappa_z$ is shown in \cref{fig:spe10_perm_z}. We can see that the contrast ratio of $\kappa^x$ is around $10^8$ and the contrast ratio of $\kappa^z$ is around $10^{11}$. We also calculate $\kappa^x\!/\kappa^z$ and present the result in \cref{fig:spe10_ratio}, which ranges from $1$ to $10^4$. The original SPE10 model contains inactive blocks, that is, cells that $\phi$ take zero. To avoid special treatments for inactive blocks, we modify the original porosity data by truncating $\phi$ at $0.05$, and the modified porosity field is displayed in \cref{fig:spe10_poro}.

The relative permeabilities $k_\mathup{w}^\mathup{r}$ and $k_\mathup{o}^\mathup{r}$ are given by the following formulas:
\[
    k_\mathup{w}^\mathup{r}=(S^*)^2,\ k_\mathup{o}^\mathup{r}=(1-S^*)^2,\ S^*=\frac{S-0.2}{0.6}.
\]
The fluid viscosities are taken as $\mu_\mathup{w}=0.3\,\mathup{cP}$ and $\mu_\mathup{o}=3.0\,\mathup{cP}$. We can hence obtain $\lambda(S)$ and $f_\mathup{w}(S)$ in \cref{eq:v,eq:trans_S}. The no-flow boundary condition is imposed. The source terms $\widetilde{q}_\mathup{w}$ and $\widetilde{q}_\mathup{o}$ are determined by well conditions. In the SPE10 model, four production wells are placed in the four corners of the $x\times y$-plane, and one injection well is located in the center of the $x\times y$-plane. All wells are vertical and penetrate completely throughout the domain. For the central injector, the injection rate is fixed as $5000\,\mathup{bbl}/\!\mathup{day}$. Therefore, for all cells that contact the central injector, we set the same $\widetilde{q}_\mathup{w}$ such that the summation equals the prescribed injection rate. For four producers, only the bottom hole pressure $p_{\scriptscriptstyle \mathup{BH}}=4000\,\mathup{psi}$ is provided. For a cell that contacts a producer, the local pressure $p_\alpha$ and source $\widetilde{q}_\alpha$ satisfy
\[
    \mathup{WI} \frac{\kappa}{\mu_\alpha}(p_{\scriptscriptstyle \mathup{BH}}-p_\alpha)=\widetilde{q}_\alpha,
\]
where $\mathup{WI}$ is called the well index, $\kappa$ is the local $\kappa^x$ or $\kappa^y$. Following the instruction in \cite{Chen2007}, we engineered the well index as
\[
    \mathup{WI} = \frac{\frac{1}{\RoundBrackets*{h^x}^2}+\frac{1}{\RoundBrackets*{h^y}^2}}{\frac{2}{\pi}\RoundBrackets*{\frac{h^y}{h^x}\log \frac{r_1}{r_\mathup{wb}}+\frac{h^x}{h^y}\log \frac{r_2}{r_\mathup{wb}}}-1},
\]
where $r_1=\sqrt{9(h^x)^2+(h^y)^2}/2$, $r_2=\sqrt{(h^x)^2+9(h^y)^2}/2$, $r_\mathup{wb}$ is termed the wellbore radius and takes the value $1.0\,\mathup{ft}$. Note now the right-hand term of \cref{eq:div_v} is coupled with $p$, which differs from the original problem \cref{eq:orgional_equation}. The initial condition in the SPE10 model is $S=0.2$, i.e., the initial water saturation equals $0.2$ in all cells.

We adopt the improved implicit pressure and explicit saturation (IMPES) method (ref.~\cite{Sheldon1959,Stone1961,Chen2007}) to solve this coupled system \cref{eq:v,eq:div_v,eq:trans_S}. More specifically, assuming that at the $n$-th time step, the saturation $S^n$ is known, we solve the system
\begin{equation}\label{eq:IMP}
    \left\{
    \begin{aligned}
        \bm{v}^{n+1}      & = -\lambda(S^n)\mathbb{K}\nabla p^{n+1}, \\
        \Div \bm{v}^{n+1} & = \widetilde{q}(p^{n+1}),
    \end{aligned}
    \right.
\end{equation}
via the discretization scheme from the velocity elimination technique to obtain the updated $p^{n+1}$ and $\bm{v}^{n+1}$. Then, set $S^{n,0}\coloneqq S^n$, for any fine element $\tau$, update the local saturation ($S^{n,m}|_\tau \rightarrow S^{n,m+1}|_\tau$) by the explicit Euler time discretization of \cref{eq:trans_S} as
\begin{equation}\label{eq:ES}
    \phi|_\tau\frac{S^{n,m+1}|_\tau-S^{n,m}|_\tau}{\delta t^{n,m}}+\frac{1}{\abs{\tau}} \int_{\partial \tau} f_\mathup{w}(S^{n,m})\bm{v}^{n+1}\cdot \bm{n} \di A=\widetilde{q}_\mathup{w}(p^{n+1}|_\tau),
\end{equation}
where $\delta t^{n,m}$ is chosen adaptively such that the maximal variance of $S^{n,m+1}$ and $S^{n,m}$ will not exceed $\mathup{DS}_{\mathup{max}}$, and $f_\mathup{w}(S^{n,m})|_{\partial \tau}$ should be determined by the upstream weighting rule (see \cite{Chen2007}). After $M$ iterations, an updated saturation field is obtained by $S^{n+1}=S^{n,M}$. As a comparison, in the classic IMPES method, the routine \cref{eq:ES} is performed only once. We realize that \cref{eq:ES} could be fully parallelized while solving \cref{eq:IMP} dominates the overall computation time.

\begin{figure}[!ht]
    \centering
    \begin{subfigure}[b]{0.49\textwidth}
        \centering
        \includegraphics[width=\textwidth]{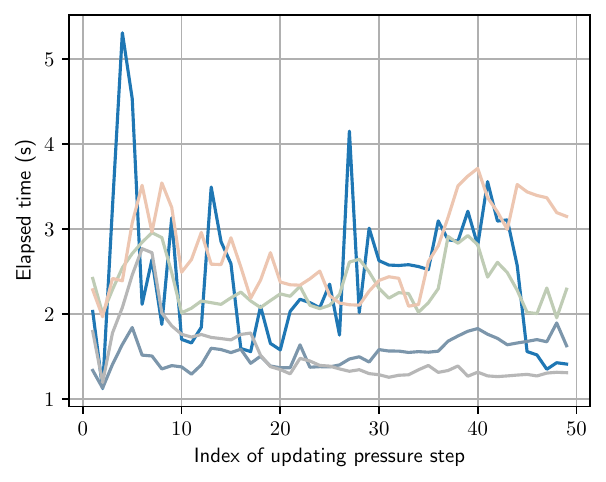}
        \caption{}\label{fig:spe10_en_time}
    \end{subfigure}
    \begin{subfigure}[b]{0.49\textwidth}
        \centering
        \includegraphics[width=\textwidth]{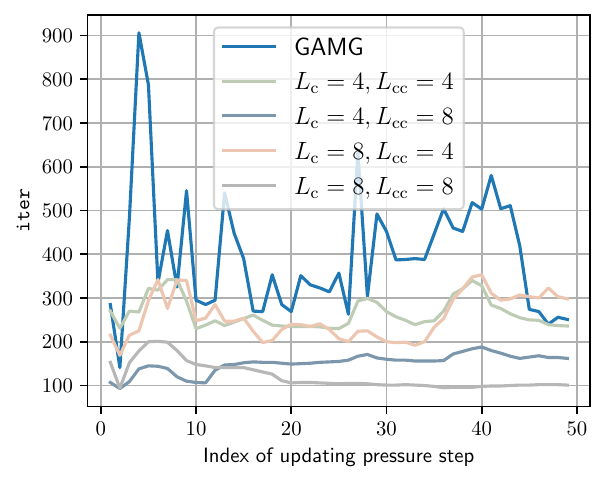}
        \caption{}\label{fig:spe10_en_iter}
    \end{subfigure}
    \begin{subfigure}[b]{0.49\textwidth}
        \centering
        \includegraphics[width=\textwidth]{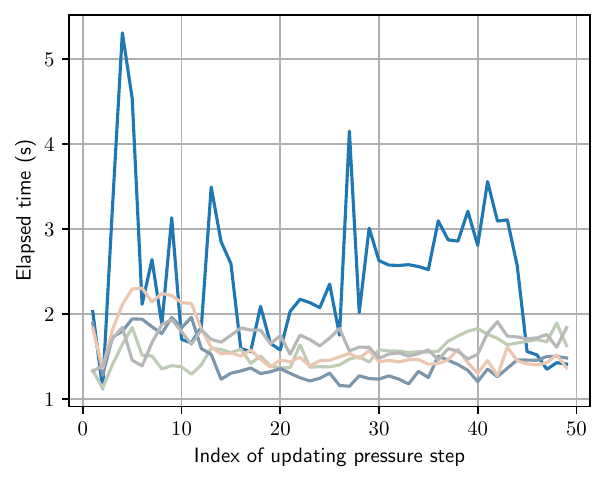}
        \caption{}\label{fig:spe10_si_time}
    \end{subfigure}
    \begin{subfigure}[b]{0.49\textwidth}
        \centering
        \includegraphics[width=\textwidth]{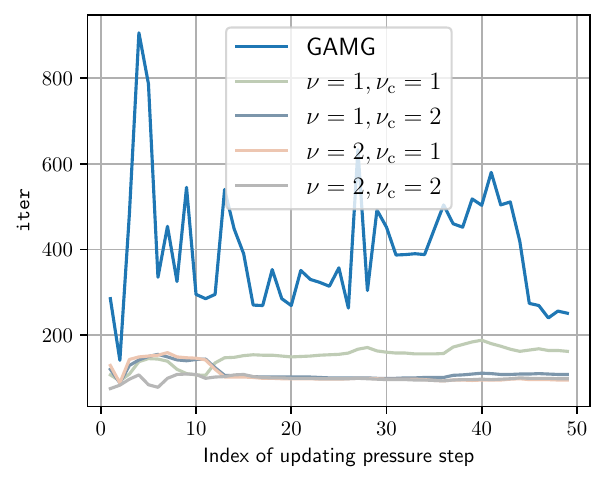}
        \caption{}\label{fig:spe10_si_iter}
    \end{subfigure}
    \caption{The elapsed time and iteration numbers of the first $50$ updating pressure steps. (a) and (b) share the same setting that $(\nu, \nu_\mathup{c})=(1,1)$ but $(L_\mathup{c}, L_\mathup{cc})$ from $\{4, 8\}\times\{4, 8\}$. (c) and (d) share the same setting that $(L_\mathup{c}, L_\mathup{cc})=(4, 8)$ but $(\nu, \nu_\mathup{c})\in \{1, 2\}\times\{1, 2\}$.}\label{fig:spe10_report}
\end{figure}

We employ the proposed spectral preconditioner in solving \cref{eq:IMP} and set the relative tolerance at $10^{-9}$. To accommodate the well condition, we modify the left-hand bilinear forms in \cref{eq:spepb,eq:spectral cc} accordingly. Except for $p^0$ the first step of pressure, solving \cref{eq:IMP} could significantly benefit from considering the previous $p^n$ as the initial guess for iterative solvers. Noting that the $\mathtt{DoF}$ of the SPE10 model equals to $60\times 220 \times 85$, we hence set $\mathtt{proc}=2 \times 7 \times 3$ to equally distribute subdomains and also avoid inter-note communications. We adopt the first strategy in constructing local coarse and coarse-coarse spaces and fix $\mathtt{sd}=4$. The rest parameters in the improved IMPES are chosen as $(\mathup{DS}_{\mathup{max}}, M)=(0.001, 50)$.  The performance is also benchmarked against the default algebraic multigrid preconditioner (GAMG) of PETSc. In the first group of experiments, we aim to investigate different choices of $(L_\mathup{c}, L_\mathup{cc})$ from $\{4, 8\}\times\{4, 8\}$ with $(\nu, \nu_\mathup{c})=(1,1)$, and the computation time and iteration number of each updating pressure step (that is, solving \cref{eq:IMP}) are shown in \cref{fig:spe10_en_time,fig:spe10_en_iter}. In the second group of experiments, we set different numbers of block Jacobi iterations within the smoothers while fixing $(L_\mathup{c}, L_\mathup{cc})$, that is, we take $(\nu, \nu_\mathup{c})\in \{1, 2\}\times\{1, 2\}$ with $(L_\mathup{c}, L_\mathup{cc})=(4, 8)$, and the results are presented in \cref{fig:spe10_si_time,fig:spe10_si_iter}. Because there is no proper initial guess available for solving $p^0$, the computing records regarding the first updating pressure step are outliers compared to other steps, and we hence exclude them in \cref{fig:spe10_en_time,fig:spe10_en_iter,fig:spe10_si_time,fig:spe10_si_iter}.

The first impression from \cref{fig:spe10_report} could be that GAMG exhibits larger variances in terms of elapsed time and iteration numbers; that is, GAMG performs less robustly compared to the proposed preconditioner. From \cref{fig:spe10_en_time,fig:spe10_en_iter}, we can see that larger $L_\mathup{c}$ and $L_\mathup{cc}$ can help to stabilize the performance, and it seems that the setting $(L_\mathup{c},L_\mathup{cc})=(4, 8)$ slightly outperforms than others regarding stability and efficiency. According to \cref{fig:spe10_si_iter,fig:spe10_si_time}, the effect of different $\nu$ and $\nu_\mathup{c}$ is rather weak, while the setting $(\nu, \nu_\mathup{c})=(1, 2)$ surpasses by a narrow margin.

\begin{figure}[!ht]
    \centering
    \begin{subfigure}[b]{0.49\textwidth}
        \centering
        \includegraphics[width=\textwidth]{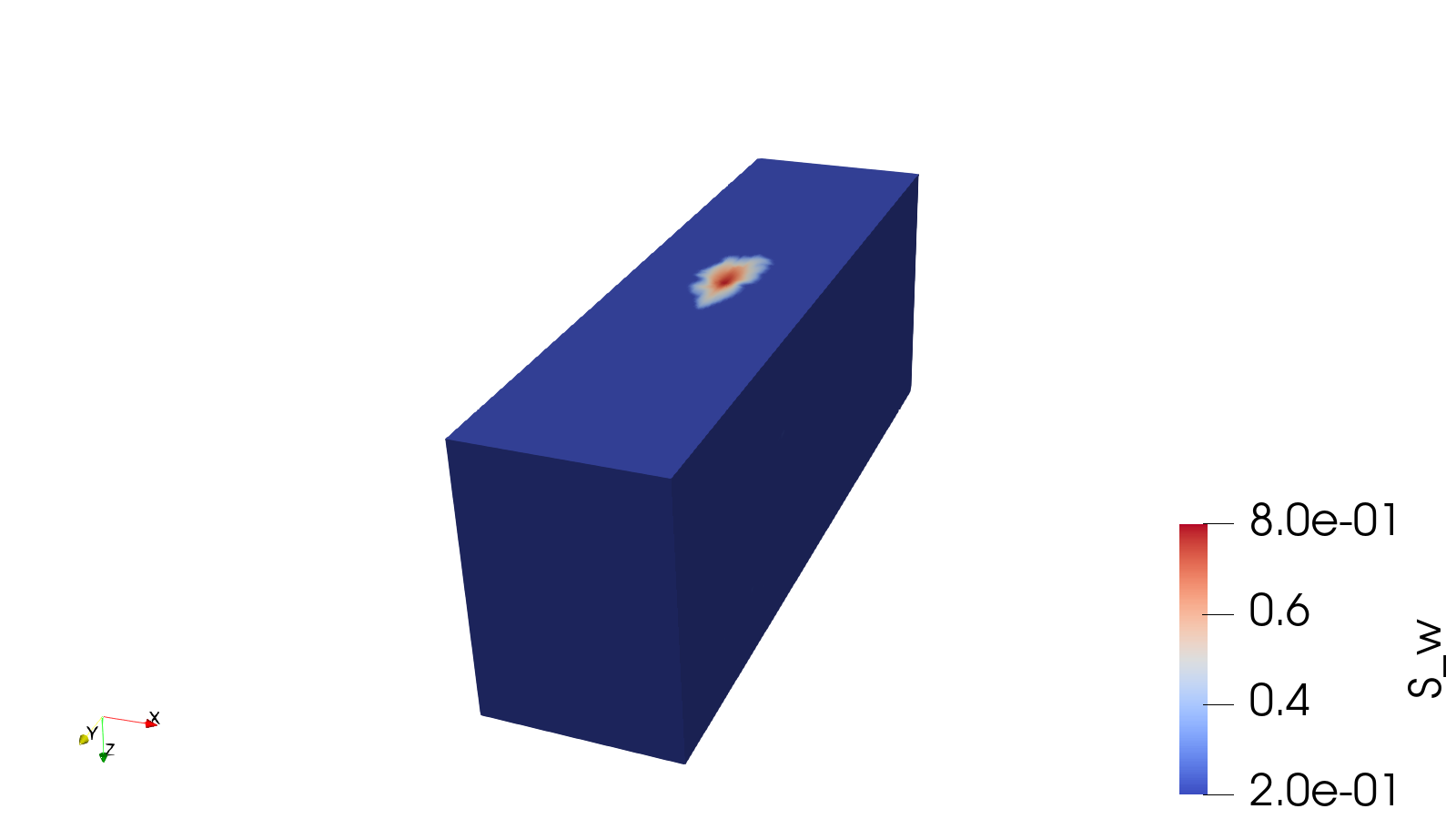}
        \caption{}
    \end{subfigure}
    \begin{subfigure}[b]{0.49\textwidth}
        \centering
        \includegraphics[width=\textwidth]{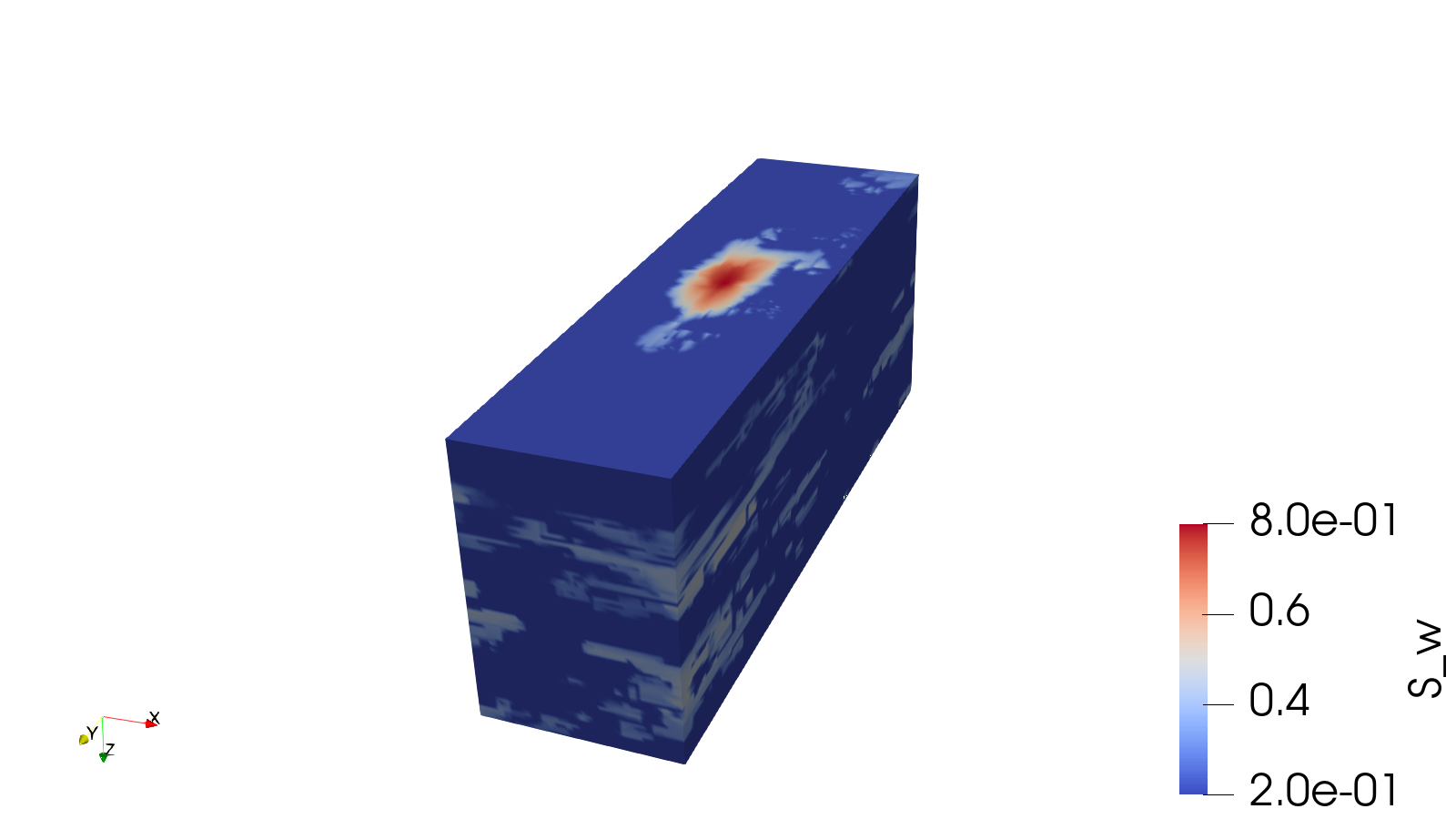}
        \caption{}
    \end{subfigure}
    \begin{subfigure}[b]{0.49\textwidth}
        \centering
        \includegraphics[width=\textwidth]{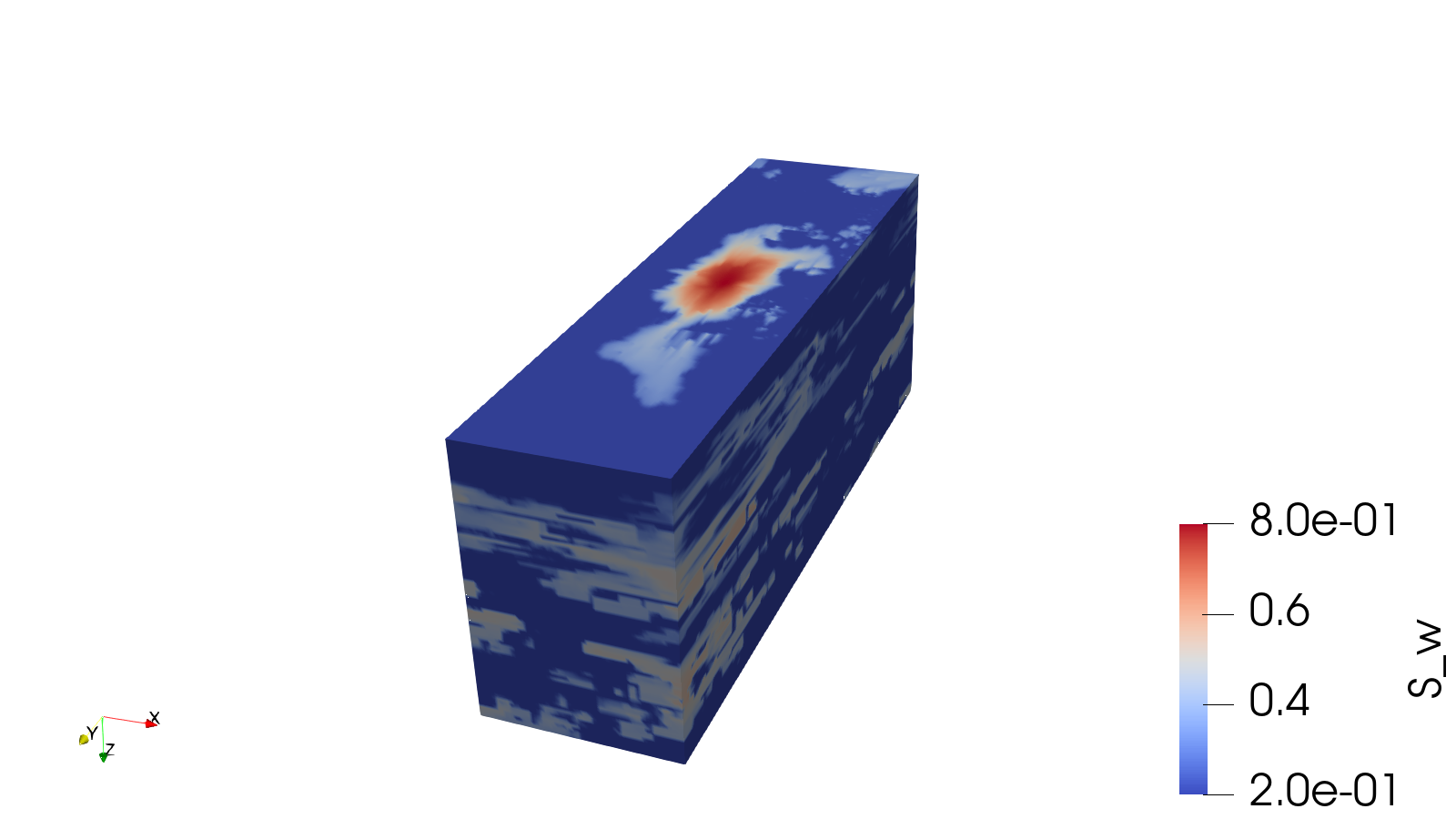}
        \caption{}
    \end{subfigure}
    \begin{subfigure}[b]{0.49\textwidth}
        \centering
        \includegraphics[width=\textwidth]{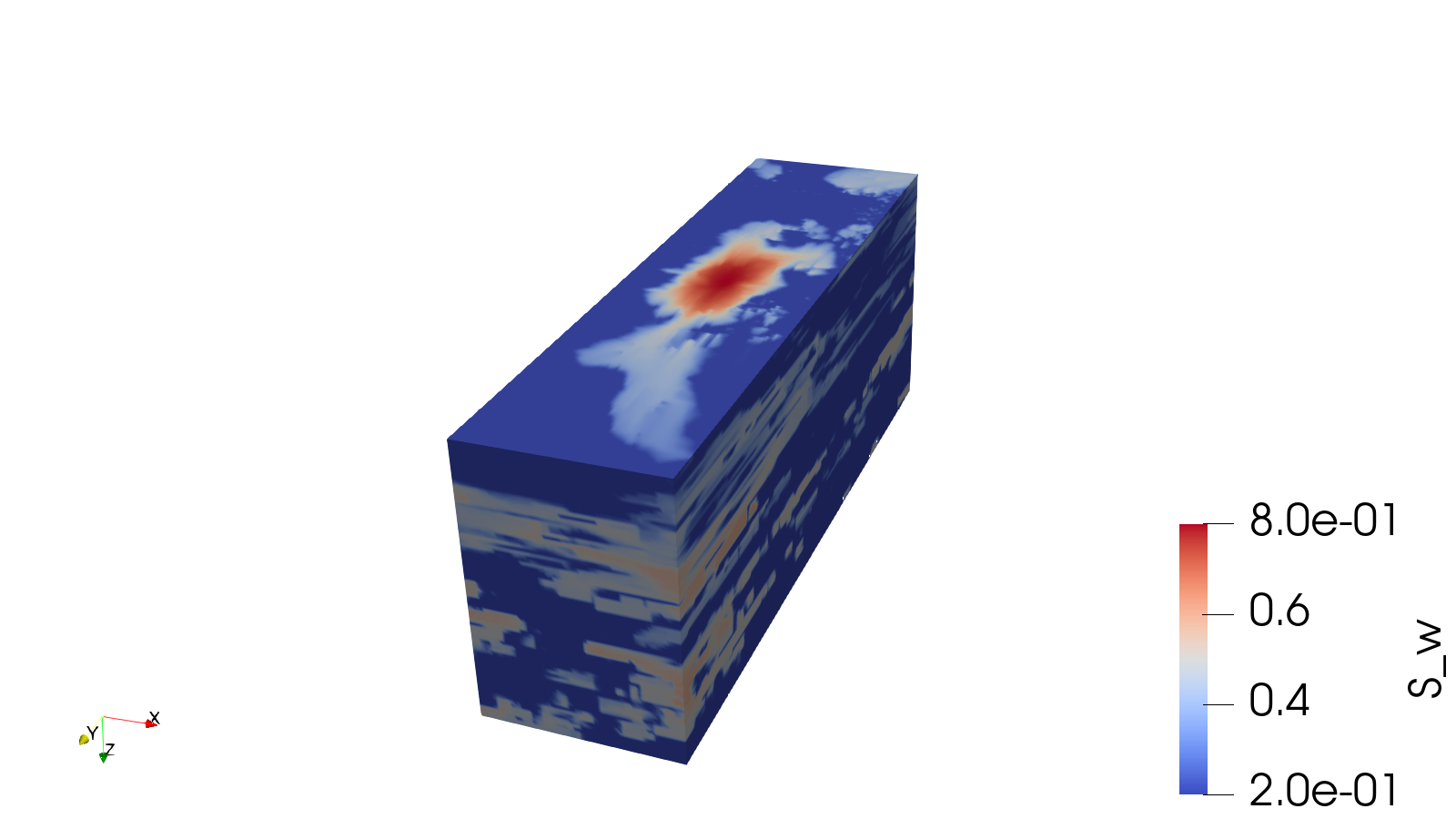}
        \caption{}
    \end{subfigure}
    \caption{The water saturation of the SPE10 model at different time: (a) $t=60.17\,\mathup{days}$, (b) $t=520.04\,\mathup{days}$, (c) $t=1211.85\,\mathup{days}$, (d) $t=1995.13\,\mathup{days}$. }\label{fig:S_w}
\end{figure}

\begin{figure}[!ht]
    \centering
    \begin{subfigure}[b]{0.49\textwidth}
        \centering
        \includegraphics[width=\textwidth]{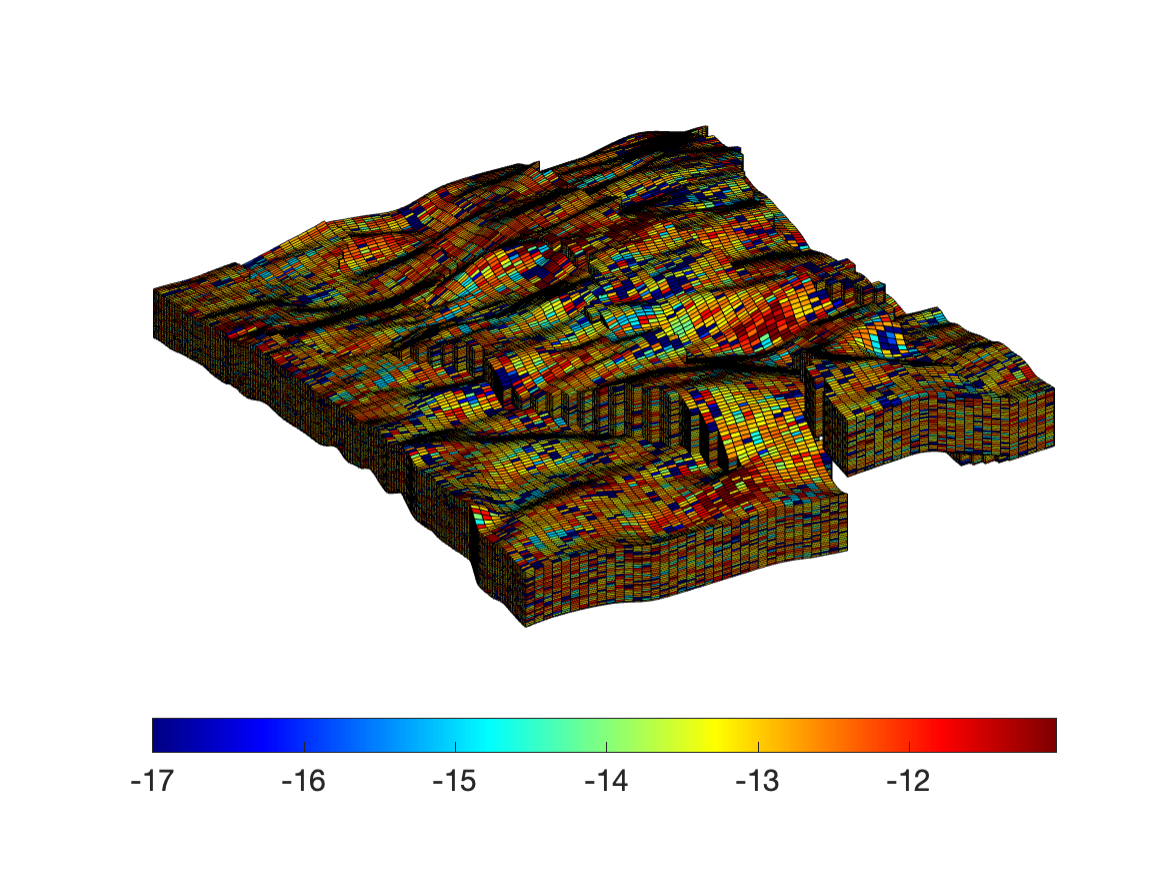}
        \caption{}
    \end{subfigure}
    \begin{subfigure}[b]{0.49\textwidth}
        \centering
        \includegraphics[width=\textwidth]{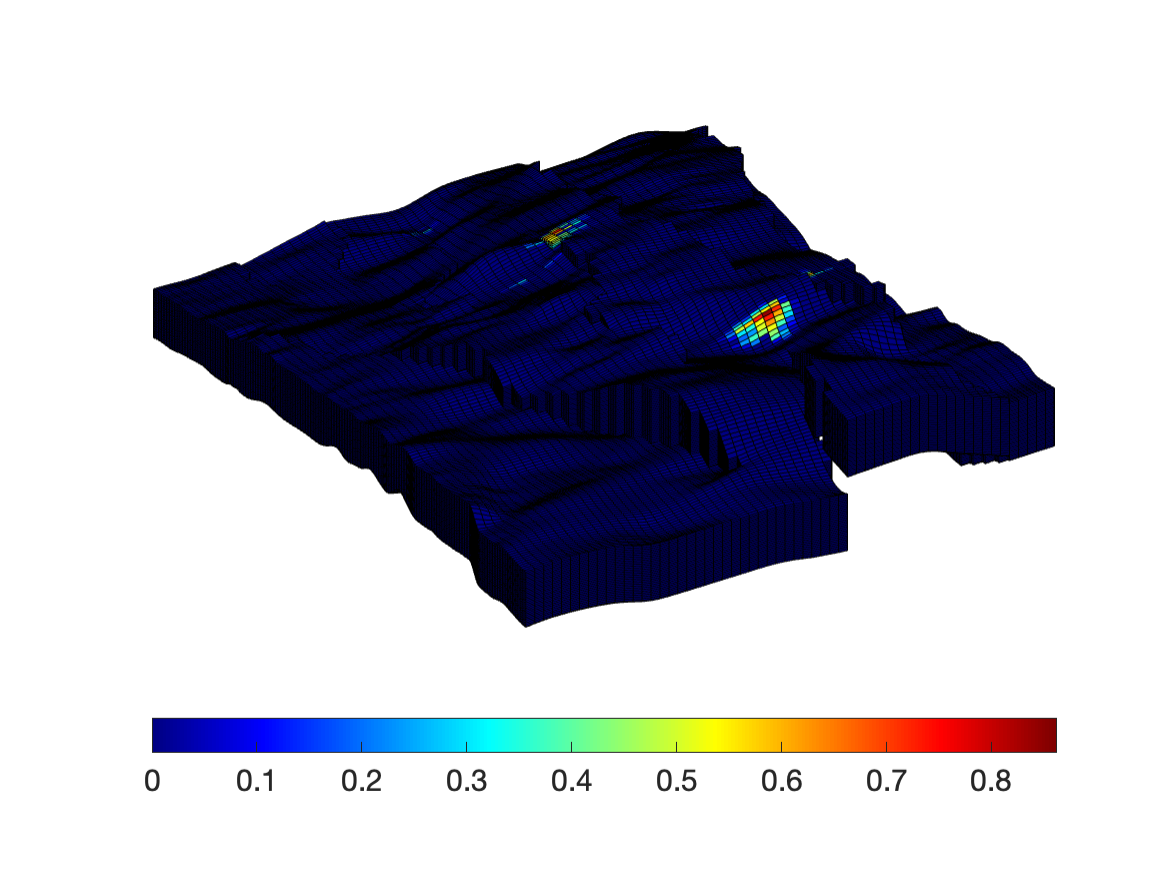}
        \caption{}
    \end{subfigure}
    \begin{subfigure}[b]{0.49\textwidth}
        \centering
        \includegraphics[width=\textwidth]{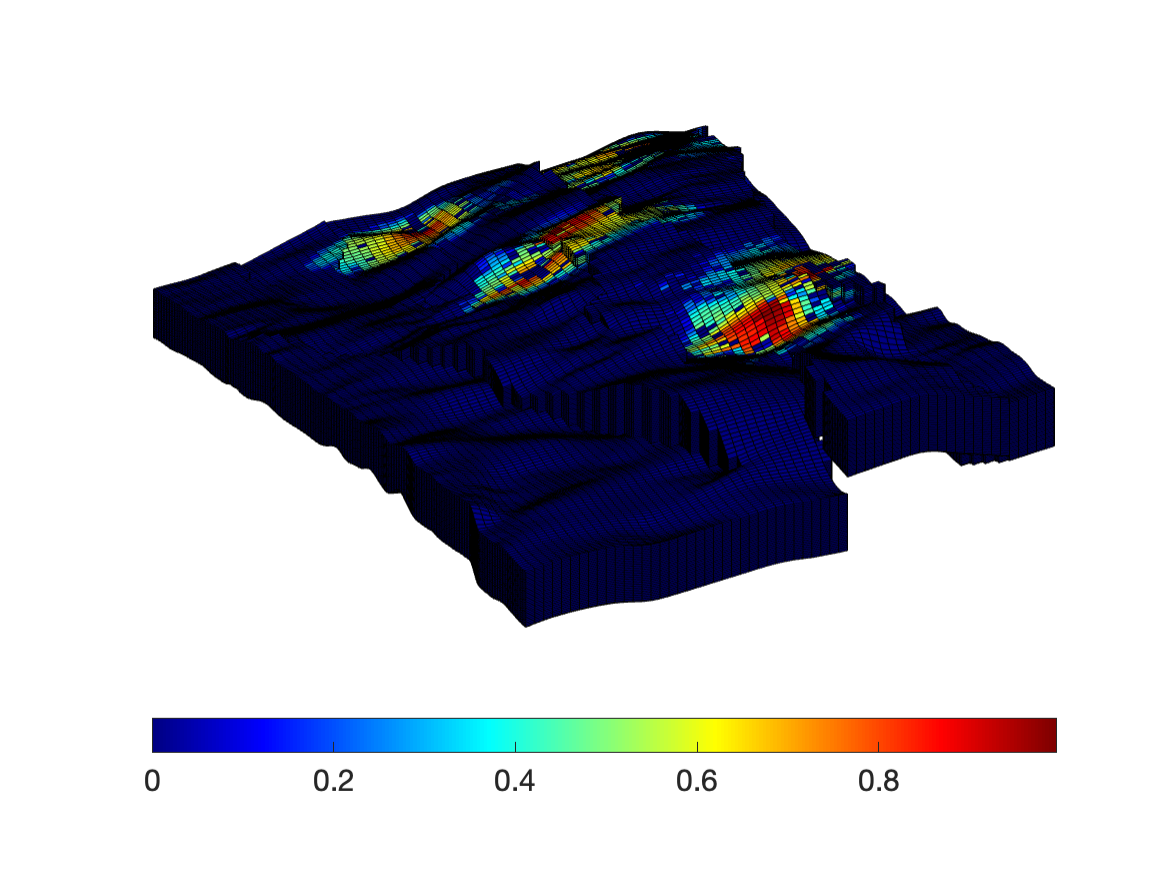}
        \caption{}
    \end{subfigure}
    \begin{subfigure}[b]{0.49\textwidth}
        \centering
        \includegraphics[width=\textwidth]{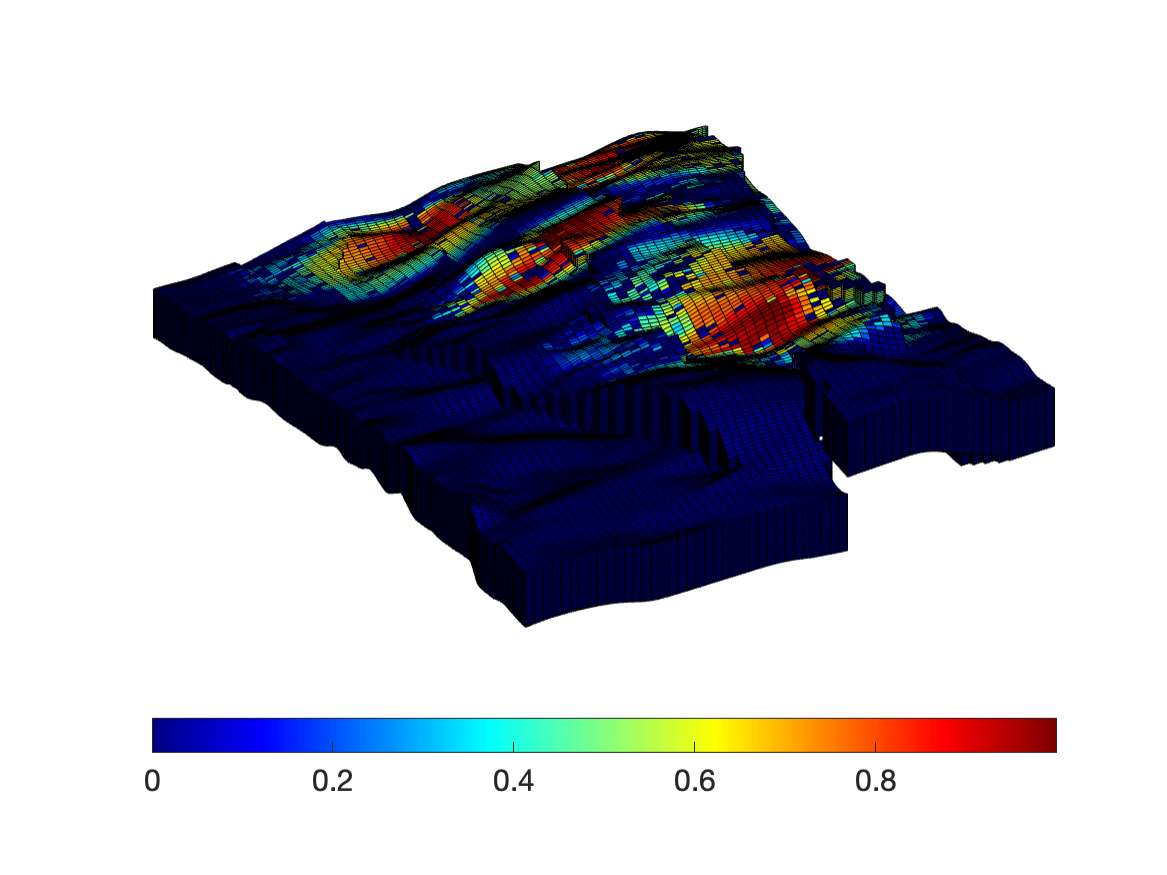}
        \caption{}
    \end{subfigure}
    \caption{
    (a) Permeability field of the Watt field model. (b)-(d) The water saturation of the Watt field model at different time: (b) $t=243.50\,\mathup{days}$, (c) $t=5113.40\,\mathup{days}$, (d) $t=29462.90\,\mathup{days}$. }\label{fig:S_w_uns}
\end{figure}

The original SPE10 project requires simulating $2000$ days of production. We are also capable of performing the $2000$-day simulation with the proposed preconditioner. The saturation of the water over several time frames is visualized in \cref{fig:S_w}. We can observe that the water saturation is much higher around the injector and that the oil is driven by water to the producers. The improved IMPES loosens the restrictions on time-step sizes, but they are not large enough for the $2000$-day simulation. 
We note with the help of \cite{lie2019introduction}, we are able to apply our method for  more realistic models such as Watt field \cite{arnold2013hierarchical}, which is based on a combination of synthetic data and real data from a North Sea oil field. 
\cref{fig:S_w_uns} show the permeability of the Watt model and saturation profiles at different time.
Currently, full implicit time methods with Constrained Pressure Residual (CPR) preconditioners (see \cite{Wallis1983,Wallis1985}) are preferred by industry solvers. 
A computing hotspot in CPR is inexactly solving the pressure equation, which the proposed preconditioner could be naturally integrated into.

\section{Conclusion}
In this paper, we developed a spectral multigrid preconditioner for Darcy flow in high-contrast media. The core component of this preconditioner is to incorporate a sequence of nested multiscale subspaces which is crucial for improving the robustness of the preconditioner. These subspaces can be constructed recursively by solving carefully designed local spectral problems. The condition number of this preconditioner is estimated and rich typical numerical examples are provided to exhibit the robustness and scalability of the algorithm. Applications for incompressible two-phase flow problems are presented. We plan to study this preconditioner for more complicated underground flow problems such as block oil simulations.

\section*{Acknowledgements}
JH's research is supported by National Key Research and Development Project of China (No. 2023YFA1011705) and National Natural Science Foundation of China (Project numbers: 12131002).
SF's research is supported by startup funding of Eastern institute of technology, Ningbo, NSFC (Project number: 12301514) and Ningbo Yongjiang Talent Programme. SF would like to thank Dr. Olav M\o yner for providing dataset of Watt model.
Research of EC is partially supported by the Hong Kong RGC General Research Fund (Project numbers: 14305222 and 14304021).




\bibliographystyle{elsarticle-num}
\bibliography{refs}
\end{document}